      \documentstyle[12pt,psfig,amscd]{amsart}
      \numberwithin{equation}{section}
\input amssym.def
\newcommand{\al}{\alpha}
\newcommand{\Ac}{{\cal A}(\phi)_{\rm conn}}
\newcommand{\C}{{\Bbb C}}
\newcommand{\cA}{{\cal A}}
\newcommand{\cAt}[1]{{\cal A}(#1)_{\rm tree}}

\newcommand{\cD}{{\cal D}}
\newcommand{\cP}{{\cal P}}
\newcommand{\cZ}{{\check Z}}
\newcommand{\D}{\Delta}
\newcommand{\e}{\varepsilon}
\newcommand{\frakg}{{\frak g}}
\newcommand{\fS}{{\cal S}}
\newcommand{\hA}{\hat{\cal A}(\phi)}
\newcommand{\hAc}{\hat{\cal A}(\phi)_{\rm conn}}
\newcommand{\hcA}{\hat{\cal A}}
\newcommand{\hoA}[1]{{\overset{\circ}{\cal A}}(#1)\hat{\ }}
\newcommand{\hD}{\hat\D}
\newcommand{\N}{{\Bbb N}}
\newcommand{\oA}{{\overset{\circ}{\cal A}}}
\newcommand{\s}{{\sigma}}
\newcommand{\Z}{{\Bbb Z}}
      \theoremstyle{plain}
      \newtheorem{thm}{Theorem}[section]
      \newtheorem{prop}[thm]{Proposition}
      \newtheorem{lem}[thm]{Lemma}
      
      \newcommand{\psdraw}[2]
         {\begin{array}{c} \hspace{-1.3mm}
         \raisebox{-4pt}{\psfig{figure=#1.ps,width=#2}}
         \hspace{-1.9mm}\end{array}}
      \theoremstyle{definition}
      \newtheorem{defn}[thm]{Definition}
      \theoremstyle{remark}
      \newtheorem{rem}[thm]{Remark}

\begin{document}

\baselineskip20pt
\title[On a universal invariant of 3-manifolds]{\ }
\author[T.\ Le., J.\ Murakami \& T.\ Ohtsuki]{\ }
\maketitle
\par\vspace{-20pt}
\renewcommand{\thefootnote}{\fnsymbol{footnote}}
\begin{center}
{\large\bf On a universal quantum invariant of 3-manifolds
\footnote{This research is supported in part
by Grand-in Aid for Scientific Research,
Ministry of Education, Science, Sports and Culture.}}
\end{center}

\begin{center}
{\sc{Thang Q.\ T.\ Le
\footnote{Department of Mathematics, SUNY at Buffalo, Buffalo, NY, USA},
Jun Murakami
\footnote{Department of Mathematics, Osaka University}
and Tomotada Ohtsuki
\footnote{Department of Mathematical and Computing Sciences,
Tokyo Institute of Technology
}}}
\end{center}

\par\bigskip

\tableofcontents

Quantum invariants of 3-manifolds was
originally proposed by Witten in \cite{Witten}.
They are given by
$$
Z_k(M,G)= \int e^{\sqrt{-1}kCS(A)}{\cal D}A
$$
which is a topological invariant of a 3-manifold $M$
with a compact Lie group $G$ and an integer $k$,
where the integral, what we call Feynman path integral,
is over all $G$ connections $A$, and
the Chern-Simons functional $CS(A)$ is given by
$$
CS(A) = \frac{1}{4\pi} \int_M Tr(A \wedge dA + \frac{2}{3} A^3).
$$
We call $Z_k(M,G)$ quantum $G$ invariant of $M$.

By perturbation theory
we expect
the asymptotic formula of $Z_k(M,G)$ for large $k$ limit,
see \cite{AxelrodSinger,AxelrodSingerII,KontsevichII}.
The formula is given by
a sum over flat connections $\rho$ as
$$
Z_k(M,G) \sim
\sum_{\rho} e^{k CS(\rho)} \frac{k^{d_{\rho}}}{R(\rho)}
\exp\left( \sum_{\Gamma}
(\frac{2\pi\sqrt{-1}}{k})^{d(\Gamma)} I_{\Gamma}(M,\rho)\right).
$$
Here $d_{\rho}$ is the dimension of
the cohomology group of the adjoint local system given by $\rho$
and $R(\rho)$ is Reidemeister torsion, and
the second sum is over all trivalent graphs $\Gamma$, and
we mean by $d(\Gamma)$ half the number of vertices of $\Gamma$, and
$I_{\Gamma}(M,\rho)$ is defined by
an integral of a product of certain $\frak g$ valued 2-forms
where the $\frak g$ bundle is twisted by $\rho$.

Here we restrict our attention to
the contribution of the trivial connection $\rho_0$
in the asymptotic formula.
Since the $\frak g$ bundle is trivial in this case,
$I_{\Gamma}(M,\rho_0)$ is decomposed into a product
of terms depending on $\frak g$ and $M$ respectively as
$$
I_{\Gamma}(M,\rho_0) = \Gamma(\frakg) Z_{\Gamma}(M)
$$
where $\Gamma(\frak g)$ is
the value obtained by \lq\lq substitute'' $\frak g$ into $\Gamma$;
for the method of substitution, see for example \cite{BarNatan}.
We note that $\Gamma(\frak g)$ depends
only on the equivalence class of $\Gamma$
with respect to the AS and IHX relations.
Further
$Z_{\Gamma}(M)$ is given by
$$
Z_{\Gamma}(M)=\int_{(x_1,\cdots,x_{2d}) \in M^{2d}- \Delta}
\prod_{\text{edges of $\Gamma$}} \omega(x_{l_i},x_{r_i})
$$
with a certain two form $\omega$ on $M \times M$,
where we put $d=d(\Gamma)$, and
the integral is over
the product of $2d$ copies of $M$
removed by its diagonal set $\Delta$, and
we take the product over all edges of $\Gamma$;
we associate the vertices of $\Gamma$ with $x_1, x_2, \cdots, x_{2d}$
and we denote by $x_{l_i}$ and $x_{r_i}$
two ends of the $i$-th edges.

On the other hand
we have a rigorous construction of
quantum $G$ invariants;
we can obtain them
by taking a sum of quantum $(\frakg, R)$ invariants of links
over a certain set of representations $R$ of $\frakg$;
for $sl_2$ case, see \cite{ReshetikhinTuraev}.
Though we fix the parameter $q$ in quantum invariants of links
to be a root of unity in the construction,
we can expect the asymptotic formula
of quantum invariants of 3-manifolds;
see \cite{Ohtsuki}
for an approach to the asymptotic formula
for quantum $SO(3)$ invariants.
For $SU(2)$ case,
relations between the formula obtained in \cite{Ohtsuki}
and the trivial connection contribution of the above perturbative expansion
are discussed in \cite{Rozansky,RozanskyII}.

Further we have
the universal Vassiliev-Kontsevich invariant
of links.
It belongs to the space of chord diagrams
consisting of solid and dashed lines
subject to the AS, STU and IHX relations.
The invariant is universal
in the sense that
it depends on neither a Lie algebra nor its representation
and we can recover
the quantum $(\frakg, R)$ invariant
by \lq\lq substituting''
the Lie algebra $\frakg$ and the representation $R$
to dashed and solid lines of the chord diagrams respectively.

Now we expect a rigorous construction of an invariant of 3-manifolds
derived from the universal Vassiliev-Kontsevich invariant,
which is universal in the sense that
it consists of dashed chord diagrams and
one should recover
the trivial connection expansion of
the perturbative expansion of quantum $G$ invariant
by substituting the Lie algebra of $G$ into the dashed lines.
In the present paper we get
an invariant of 3-manifold
from the universal Vassiliev-Kontsevich invariant,
which is an infinite linear sum of chord diagrams,
that is,
it belongs to
the space of chord diagrams
consisting of dashed lines (i.e. trivalent graphs)
subject to the AS and IHX relations.
The space is graded by the number of vertices.
The simplest one is the theta curve
and we show that the coefficient of it corresponds to the Casson invariant.

In Section 1 we review a definition of
the universal Vassiliev-Kontsevich invariant
and show properties of it.
In Section 2 we give a map to remove solid lines,
and we define a universal quantum invariant using the map
in Sections 3 and 4.
We show some properties of the invariant in Section 5.

In the previous papers \cite{LMMO,LMMOII}
only the coefficient of the theta curve was obtained
from the universal Vassiliev-Kontsevich invariant.
The results of the present paper were announced in \cite{LMO}.

The authors would like to thank the warm hospitality of
Matematisk Institut, \AA{rhus} Universitet
for our stay in the summer of 1995,
where a part of this work was done.
They are also grateful to Hitoshi Murakami for valuable conversations
and encouragement.

\section{Modified universal Vassiliev-Kontsevich invariant}

We will review
a construction of the universal Vassiliev-Kontsevich invariant
of framed oriented links,
and the definition of modified one given in \cite{LeMurakamiII}.
The modified invariant is well behaved
under Kirby move II.
We also show other properties of the invariant in this section.

\subsection{Chord diagrams}
A {\it uni-trivalent graph}
is a graph every vertex of which is either univalent
or trivalent. A uni-trivalent graph is {\it oriented} if at each trivalent
 vertex a cyclic order of edges is fixed.

Let $X$ be a compact oriented 1-dimensional manifold
whose components are labeled.
A {\it chord diagram} with support $X$ is the manifold $X$ together with an
oriented uni-trivalent graph
 whose univalent vertices are on $X$; and the graph does not have any connected
component
homeomorphic to a circle.
Note that our definition of a chord diagram is more general than
that of \cite{BarNatan,LeMurakami}.
In Figures components of $X$ are depicted by solid lines,
while the graph is depicted by dashed lines.
There may be connected components of the  dashed graph which do
not have univalent vertices.
Each chord diagram has a natural topology.
Two chord diagrams $D,D'$ on $X$ are regarded
as equal if there is a homeomorphism
$f:D\to D'$ such that $f|_X$ is a homeomorphism
of $X$ which preserves components and orientation.

\begin{figure}[htpb]
$$ \psdraw{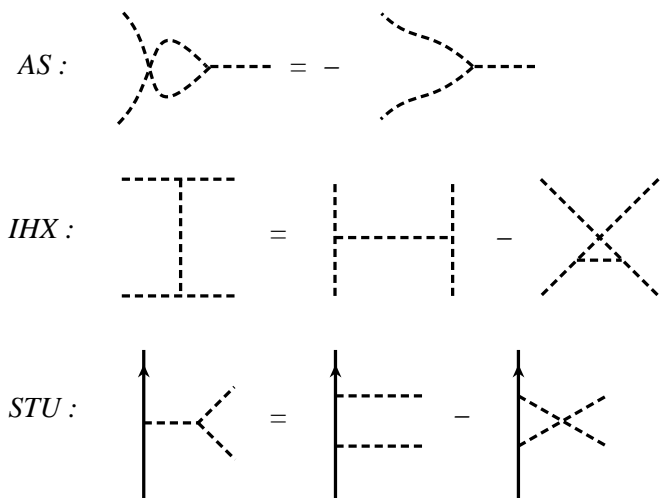}{3.6in} $$
\caption{The AS, IHX and STU relations}\label{fig.ASIHXSTU}
\end{figure}

Let $\cA(X)$ be the vector space over $\C$
spanned by chord diagrams with support $X$,
subject to the AS, IHX and STU relations
shown in Figure \ref{fig.ASIHXSTU}.
The {\it degree} of a chord diagram is
half the number of vertices of the dashed graph.
Since the relations AS, IHX and STU respect the degree,
there is a grading on $\cA(X)$ induced by this degree.
We denote by $\hcA(X)$ the completion of $\cA(X)$
with respect to the degree.
Note that
each map on $\cA(X)$ defined in the following
has the natural extension of it on $\hcA(X)$,
and that
we use the same notation for the corresponding maps
on $\cA(X)$ and $\hcA(X)$.

Suppose $C$ is a component of $X$.
Reversing the orientation of $C$, from $X$  we
get $X'$. Let $S_{(C)}: \cA(X)\to\cA(X')$ be the linear mapping which
transfers every chord diagram $D$ in $\cA(X)$ to $S_{(C)}(D)$ obtained
from $D$ by reversing the orientation of $C$ and multiplying  by
$(-1)^m$, where $m$ is the number of vertices of the dashed graph
on the component $C$.

Replacing $C$ by $2$ copies of $C$,
from $X$ we get $X^{(2,C)}$,
with a projection $p:X^{(2,C)}\to X$.
If $x$ is a point on $C$ then $p^{-1}(x)$ consists of $2$ points,
while if $x$ is a point of other components,
then $p^{-1}(x)$ consists of one point.
Let $D$ be a chord diagram on $X$,
with the dashed graph $G$.
Suppose that there are $m$ univalent vertices of $G$ on $C$.
Consider all possible new chord diagrams on $X^{(2,C)}$
with the same dashed graph $G$ such that
if a univalent vertex of $G$ is attached to a point $x$ on $X$ in $D$,
then this vertex is attached to a point in $p^{-1}(x)$
in the new chord diagram.
There are $2^m$ such chord diagrams,
and their sum is denoted by $\Delta_{C}(D)$.
It is easy to check that these linear mappings $S_{(C)}$ and $\Delta_{(C)}$
are well-defined on $\cA(X)$,
and naturally extended to the maps on $\hcA(X)$.

Suppose that $X$ and $X'$ have distinguished components
$C$ and $C'$ respectively,
and that $X$ consists of loop components only.
Let $D\in\cA(X)$ and $D'\in\cA(X')$ be two chord diagrams.
{}From each  of $C$ and $C'$ we remove a small arc
which does not contain any vertices.
The remaining part of $C$ is an arc which we glue  to $C'$ in
the place of the removed arc such that the orientations are compatible.
The new chord diagram is called the
{\it connected sum of $D$ and $D'$ along the distinguished components};
it does not depend on the
locations of the removed arcs, which follows from the STU relation and
the fact that all components of $X$ are loops. The proof is
the same as in the case $X=X'=S^1$ as in \cite{BarNatan}.

We define a co-multiplication $\hat \Delta$
in $\cA(X)$ and $\hcA(X)$ as follows.
A {\it  chord sub-diagram} of a chord diagram $D$ with dashed graph $G$
 is any chord
 diagram obtained from $D$ by removing
some  connected components of $G$. The {\it complement chord sub-diagram}
of  a chord sub-diagram
$D'$
is  the chord  sub-diagram obtained by removing components of $G$
which are in $D'$.
We define
$$\hat \Delta (D)=\sum D'\otimes D''.$$
Here the sum is over all chord sub-diagrams $D'$ of $D$, and $D''$ is
 the complement
of $D'$. This co-multiplication is co-commutative.

\subsection{Associator}
Let $\C<<A,B>>$ be the algebra over $\C$ of
 all formal power series in two non-commutative symbols $A,B$.
We are going to define an element $\varphi\in \C<<A,B>>$,
known as the Drinfeld associator.

We put
\begin{equation}
\zeta(i_1,\ldots,i_k)=\sum_{n_1<\ldots<n_k\in\N}\frac{1}
{n_1^{i_1}\cdots n_k^{i_k}},
\end{equation}
for natural numbers $i_1,\ldots,i_k$ satisfying $i_k\ge2$.
These values, called multiple zeta values, have recently gained
 much attention among number theorists.
In what follows
bold letters ${\bold p}, {\bold q}, {\bold r}, {\bold s}$
stand for non-negative multi-indices.
For a multi-index ${\mathbf{p}}=(p_1,\cdots,p_k)$,
we call $k$ the {\it length} of $\mathbf{p}$.
Let ${\mathbf{1}}_k$ be the multi-index consisting of $k$ letters 1.
We denote $\sum p_i$ by $|{\mathbf{p}}|$.
For two multi-indices $\mathbf{p}$ and $\mathbf{q}$ of the same length $k$,
we put
$\eta ({\mathbf{p}};{\mathbf{q}})=0$ if one of $p_i,q_i$ is 0,
$$
\eta ({\mathbf{p}};{\mathbf{q}})
\ =\ \zeta ({\mathbf{1}}_{p_1-1},q_1+1,{\mathbf{1}}_{p_2-1},
q_2+1,\dots ,{\mathbf{1}}_{p_k-1},q_k+1),
$$
otherwise.
Further we set two notations by
\begin{align*}
(A,B)^{({\mathbf{p}},{\mathbf{q}})}
&=A^{p_1}B^{q_1}A^{p_2}B^{q_2}\dots A^{p_k}B^{q_k}, \\
\binom{\mathbf{p}}{\mathbf{q}}
&=\binom{p_1}{q_1}\binom{p_2}{q_2}\cdots \binom{p_k}{q_k}.
\end{align*}
Using the above notations we define the {\it associator} $\varphi$ by
\[
\varphi(A,B)=1+\sum_{\mathbf{p},\mathbf{q},\mathbf{r},\mathbf{s}\ge 0}
(-1)^{|{\bold{r}}|+|{\bold{q}}|}\eta ({\mathbf{p}}+{\mathbf{r}};
{\mathbf{q}}+{\mathbf{s}})
\,\binom{\mathbf{p+r}}{\mathbf{r}}
\binom{\mathbf{q+s}}{\mathbf{s}}B^{|{\mathbf{s}}|}
(A,B)^{({\mathbf{p}},{\mathbf{q}})}\,A^{|{\mathbf{r}}|}.
\]
Here the sum is over all multi-indices $\mathbf{p,q,r,s}$
of the same length $k$ where $k=1,2,3,\cdots$.
Note that there exists the inverse of $\varphi(A,B)$;
we denote it by $\varphi^{-1}(A,B)$.

\subsection{The universal Vassiliev-Kontsevich invariant
for framed oriented links}

In this subsection,
we will define the universal Vassiliev-Kontsevich invariant
$\hat Z(L)$ of a framed oriented link $L$.
Before defining it,
we define an invariant $Z(\cD)$ of a link diagram $\cD$.

\begin{figure}[htpb]
$$ \psdraw{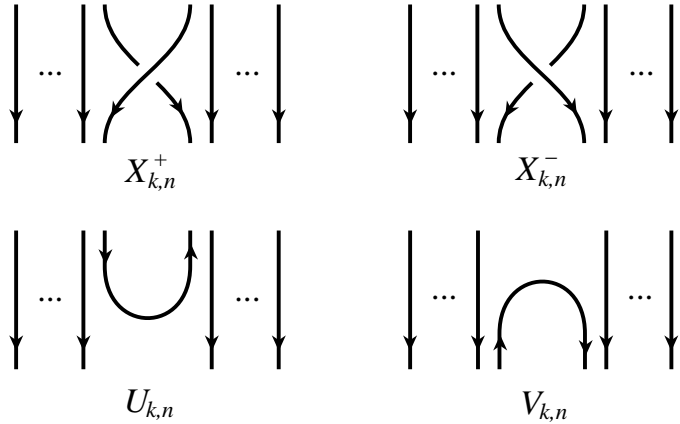}{3.6in} $$
\caption{Elementary tangle diagrams}\label{elementary}
\end{figure}

Suppose $\cD$ is an {oriented $l$-component link diagram} in a plane
with a fixed coordinate system.
Using horizontal lines one can decompose $\cD$
into elementary tangle diagrams.
Here we mean an {\it elementary tangle diagram}
is one of the tangle diagrams shown in Figure \ref{elementary},
maybe with reverse orientation on some components.
We number the components of the elementary tangle diagrams from
left to right.
There are $n$ components in $X^{\pm}_{k,n}$, and the crossing is
 on the $k$-th and $(k+1)$-th components. There are $n-1$ components
in $U_{k,n}, V_{k,n}$, and the non-straight component is numbered by $k$.
We define $Z(\cD)\in\hcA(\coprod^lS^1)$ as follows.

\begin{figure}[htpb]
$$ \psdraw{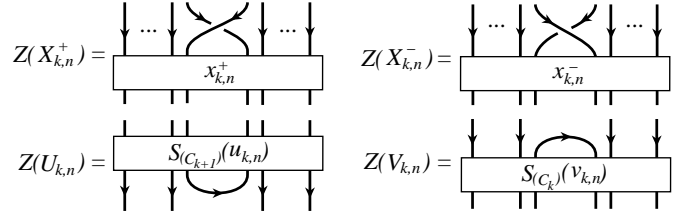}{3.6in} $$
\caption{The definition of $Z(T)$
for each elementary tangle diagram $T$}\label{fig.definexuv}
\end{figure}

Firstly we define $Z(T)$ for each elementary tangle diagram $T$.
When $X$ is $n$ vertical straight numbered lines
with downward orientations,
$\hcA(X)$ is denoted by $\cP_n$.
All the $\cP_n$ are algebras:
the product of two chord diagrams $D_1$ and $D_2$
is obtained by placing $D_1$ on the top of $D_2$.
The algebra $\cP_1$ is commutative (see \cite{BarNatan,Kontsevich}).
Consider the following  element $x_{k,n}$ in $\cP_n$
$$
x_{k,n}^+=\varphi^{-1}(\frac{1}{2\pi\sqrt{-1}}\sum_{i=1}^{k-1}\Omega_{i,k},
\frac{\Omega_{k,k+1}}{2\pi\sqrt{-1}})
\exp(\frac{\Omega_{k,k+1}}{2})\varphi(\frac{1}{2\pi\sqrt{-1}}\sum_{i=1}^{k-1}
\Omega_{i,k+1},
\frac{\Omega_{k,k+1}}{2\pi\sqrt{-1}}).
$$
Here $\varphi$ is defined in the previous subsection, and $\Omega_{i,j}$
is the chord diagram in $\cP_n$ with the dashed graph being a line
connecting the $i$-th and the $j$-th string.
Similarly we put
\begin{align*}
x_{k,n}^-
&=\varphi^{-1}(\frac{1}{2\pi\sqrt{-1}}\sum_{i=1}^{k-1}\Omega_{i,k},
\frac{\Omega_{k,k+1}}{2\pi\sqrt{-1}})
\exp(-\frac{\Omega_{k,k+1}}{2})\varphi(\frac{1}{2\pi\sqrt{-1}}
\sum_{i=1}^{k-1}\Omega_{i,k+1},
\frac{\Omega_{k,k+1}}{2\pi\sqrt{-1}}), \\ 
u_{k,n}
&=\varphi^{-1}(\frac{1}{2\pi\sqrt{-1}}\sum_{i=1}^{k-1}\Omega_{i,k},
\frac{\Omega_{k,k+1}}{2\pi\sqrt{-1}}), \\
v_{k,n}
&=\varphi(\frac{1}{2\pi\sqrt{-1}}\sum_{i=1}^{k-1}\Omega_{i,k},
\frac{\Omega_{k,k+1}}{2\pi\sqrt{-1}}).
\end{align*}
We define $Z(X_{k,n}^+)$ as the element in $\hcA(X_{k,n}^+)$
obtained by placing $X^+_{k,n}$ without any dashed graph
on the top of $x_{k,n}^+$.
Here we also use the notation $X_{k,n}^+$
for the set of solid lines of $X_{k,n}^+$.
Similarly,
$Z(X_{k,n}^-)$ is the element in $\hcA(X_{k,n}^-)$
obtained by placing $X^-_{k,n}$ on the top of $x_{k,n}^-$.
Further $Z(U_{k,n})$ is obtained by placing $U_{k,n}$
on the bottom of $S_{C_{k+1}}(u_{k,n})$,
where $C_{k+1}$ is the $(k+1)$-th strings
of the support of chord diagrams in $\cP_n$,
and $S_{C_k}$ is defined in the previous subsection.
Furthermore $Z(V_{k,n})$ is obtained by placing $V_{k,n}$
on the top of $S_{C_k}(v_{k,n})$.
We also show pictures of the definition
in Figure \ref{fig.definexuv}.

Secondly we define $Z(T)$
for each elementary tangle diagram $T$ with arbitrary orientation
by using $Z(T)=S_C(Z(T'))$ some times,
where $T$ is obtained from $T'$
by reversing the orientation of a component $C$.

Lastly we define $Z(\cD) \in \hcA(\coprod^l S^1)$
for an oriented $l$-component link diagram $\cD$.
We can decompose $\cD$ into
elementary tangle diagrams $T_1,T_2,\dots,T_m$
by horizontal lines,
counting from top to bottom.
We set
$Z(\cD)=Z(T_1)\times Z(T_2)\times \dots\times Z(T_m)$.
Here, for chord diagrams $D_i$ composing $Z(T_i)$ ($i=1,2,\dots,m$),
we mean by $D_1\times\dots\times D_m$
the chord diagram on $\coprod^l S^1$
obtained by placing $D_1$ on the top of $D_2$,
and placing the union on top of $D_3$, and so on.
Note that the supports of
$D_1,D_2,\dots,D_m$ can be glued together,
and the result consists of $l$ solid loops.
It is known that $Z(\cD)$ is defined uniquely
for an oriented link diagram $\cD$,
not depending on the decomposition of $\cD$,
i.e. more precisely,
not depending on an isotopy of the link diagram $\cD$ in the plane
which preserves the maximal points of $\cD$.

\begin{figure}[htpb]
$$ \psdraw{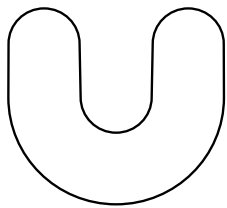}{1.2in} $$
\caption{The link diagram $U$}\label{U}
\end{figure}

Suppose now that $L$ is a framed oriented link,
represented by a link diagram $\cD$
with blackboard framing.
Further suppose that the $i$-th components of $\cD$ has
$s_i$ maximal points
with respect to the height function of the plane.
We define an invariant $\hat Z(L)$ by
$$
\hat Z(L)=Z(\cD)\#(\nu^{s_1}\otimes\dots\otimes \nu^{s_l})
\in \hcA(\coprod^l S^1),
$$
where
we put $\nu=Z(U)^{-1}$ for the link diagram $U$ shown in Figure \ref{U},
and we mean by the above formula that
$\hat Z(L)$ is obtained from $Z(\cD)$ by successively taking connected
sum with $\nu^{s_i}$ along the $i$-th component.
(Note that
a notation $\hat Z_f(L)$ was used
instead of $\hat Z(L)$
in the previous papers \cite{LeMurakami,LeMurakamiII}.)
The invariant $\hat Z(L)$ of a framed link $L$ is well-defined,
not depending on the choice of its link diagram $\cD$,
see \cite{LeMurakami}.
Note that,
for the trivial knot $K$ with framing 0,
$\hat Z(K)$ is not trivial;
in fact, we have $\hat Z(K)=\nu$.

The invariant $\hat Z(L)$ is an invariant of framed oriented links
such that
it contains in itself all Vassiliev invariants of framed oriented links,
and it is a generalization of the Kontsevich integral.
We call this invariant
{\it the universal Vassiliev-Kontsevich invariant}.

\subsection{Parallel of framed links}
The following was proved in \cite{LeMurakamiII}.

\begin{prop}\label{prop.ZandDcommute}
Let $C$ be a component of an oriented framed link $L$.
\begin{enumerate}
\item
Let $L^{(2,C)}$ be the link obtained
from $L$ by replacing $C$ by two push-offs of $C$ using the frame.
Then we have the following formula;
\begin{equation}\label{eq.parallel}
\hat Z(L^{(2,C)})=\Delta_{(C)}(\hat Z(L)).
\end{equation}
\item
Let $L'$ be obtained from $L$
by reversing the orientation of $C$.
Then we have the following formula;
\begin{equation}\label{eq.revori}
\hat Z(L')=S_{(C)}(\hat Z(L)).
\end{equation}
\end{enumerate}
\end{prop}

Now we put
$$
\check Z(L)=\hat Z(L)\# (\nu\otimes\dots\otimes \nu)
\in \hcA(\coprod^l S^1).
$$
This formula means that
$\check Z(L)$ is obtained from $\hat Z(L)$ by successively
taking connected sum with $\nu$ along every component of $L$.
It is easy to see
that Proposition \ref{prop.ZandDcommute}
is also valid if we replace $\hat Z$ by $\cZ$
in the statement of the proposition.

Suppose that $X$ consists of $n$ components $C_1,\dots,C_n$.
We denote by $X\sqcup X$ the disjoint union of two copies of $X$.
We define a linear map:
$$
p : {\hcA}(X\sqcup X) \to {\hcA}(X)^{\otimes 2}
$$
as follows.
If $D$ is a
chord diagram having dashed graph connecting the two copies  of $X$,
then we put $p(D)=0$.
Otherwise, $D$ splits into a disjoint union
of two chord diagrams $D_1$ and $D_2$
on the first and the second copies of $X$ respectively,
and we put $p(D)=D_1 \otimes D_2$.

Then by definition we have
\begin{equation}\label{co-mul}
p \circ \Delta_{(C_1,\dots,C_n)}(D)
=\hat{\Delta}(D),   \label{x1}
\end{equation}
where $\Delta_{(C_1,\dots,C_n)}:\hcA(X)\to\hcA(X\sqcup X)$ is
the mapping obtained by successively applying
$\Delta_{(C_i)},i=1,\dots,n$.

\begin{thm}\label{thm.bunshin}
Let $L$ be an oriented framed link.
Then the following formula holds;
$$
\hD(\cZ(L))=\cZ(L) \otimes \cZ(L).
$$
\end{thm}

\begin{pf}
By the equations (\ref{eq.parallel}) and (\ref{x1}),
the left hand side of the required formula is equal to
$p(\cZ(L^{(2)}))$.
Identifying $\hcA(X) \otimes \hcA(X)$
with a subset of $\hcA(X \sqcup X)$ naturally,
the right hand side of the required formula becomes equal to
$\cZ(L \circ L)$,
where $L \circ L$ is the split union of two copies of $L$,
i.e., the disjoint union not winding with each other.

Hence it is sufficient to show that
$p(\cZ(L^{(2)}))$ is equal to $\cZ(L \circ L)$.
Note that we can obtain $L \circ L$ from $L^{(2)}$
by taking crossing changes
between the first and the second copies of $L$.
In this procedure $\cZ(L^{(2)})$ changes to $\cZ(L \circ L)$
only in terms which have
dashed graphs connecting the two copies of $X$,
since the difference of
a crossing change is locally given by
\begin{align*}
x_{k,n}^+ - x_{k,n}^-
&=
\varphi^{-1}(\ldots)
\left(
\exp(\frac{\Omega_{k,k+1}}{2}) - \exp(-\frac{\Omega_{k,k+1}}{2})
\right)
\varphi(\ldots)
\\
&=
\varphi^{-1}(\ldots)
\left(
\Omega_{k,k+1} + \text{higher terms with respect to $\Omega_{k,k+1}$}
\right)
\varphi(\ldots).
\end{align*}
Therefore, after the procedure,
$\cZ(L^{(2)})$ becomes an element of $\hcA(X \sqcup X)$
which is different from $\cZ(L^{(2)})$
only in terms with dashed graphs connecting the two copies of $X$.
Further the element should vanish in such terms,
since it is equal to the invariant of a split union.
Hence it must be equal to $p(\cZ(L^{(2)}))$;
it is also equal to $\cZ(L \circ L)$.
\end{pf}

\subsection{The change under Kirby move II}

\begin{figure}[htpb]
$$ \psdraw{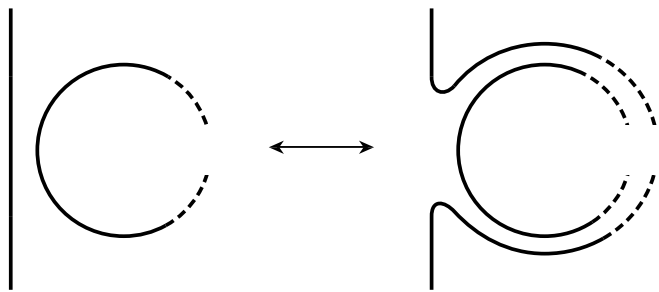}{3.6in} $$
\caption{Kirby move II (the handle slide move)}\label{fig.handleslide}
\end{figure}

We have the following proposition proved in \cite{LeMurakamiII}.
By this proposition we can get the change
of $\cZ(L)$ under Kirby move II.
Throughout the present paper
we mean by Kirby move II
the handle slide move defined by Kirby in \cite{Kirby};
we show a picture of the move in Figure \ref{fig.handleslide}.

\begin{figure}[htpb]
$$ \psdraw{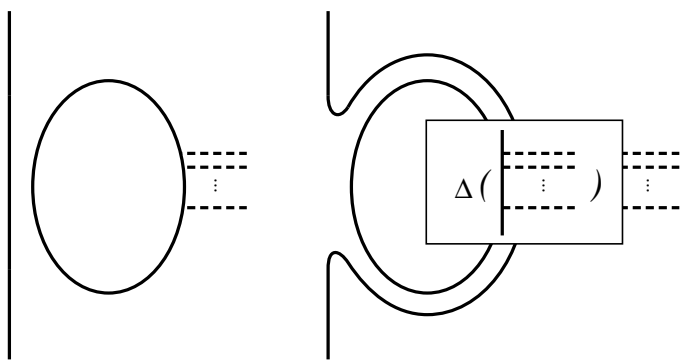}{3.6in} $$
\caption{The change under Kirby move II}\label{fig.KII}
\end{figure}

\begin{prop}[\cite{LeMurakamiII}]\label{prop.KII}
Let $L$ be an oriented framed link,
and $L'$ a framed link obtained from $L$
by Kirby move II which preserves the orientation.
Then $\cZ(L')$ can be obtained from $\cZ(L)$
by replacing the left picture in Figure \ref{fig.KII}
with the right picture.
\end{prop}

\section{Replacing solid circles by dashed lines}

Our aim of this section
is to construct the series of the maps $\iota_n$,
whose definition is given in the end of this section.
The universal Vassiliev-Kontsevich invariant
belongs to the space
consisting of both of solid and dashed lines.
When we consider quantum $({\frak g}, R)$ invariants of links
for a Lie algebra $\frak g$ and its representation $R$,
the solid and dashed lines corresponds to
$R$ and $\frak g$ respectively.
Further, we know that
no particular representation $R$ is specified
in quantum $\frak g$ invariants of 3-manifolds.
It implies that
our obstruction in constructing
invariants of 3-manifolds
from the universal Vassiliev-Kontsevich invariant
might be the existence of solid lines.
We will remove solid lines
from the space of chord diagrams
by the maps $\iota_n$
when we construct invariants of 3-manifolds
in the following section.

\subsection{Chord diagrams with finite support}

Let $X$ be a finite set
consisting of $m$ ordered points named $0,1,2,$ $\cdots,m-1$.
We consider chord diagrams with support $X$,
that is, oriented uni-trivalent graphs
whose $m$ univalent vertices are on the $m$ fixed points respectively.
We denote by $\cA(m)$
the vector space over $\C$ spanned by such chord diagrams
subject to the AS and IHX relations.
Further we denote by $\cAt{m}$ the vector subspace of $\cA(m)$
spanned by connected and simply connected graphs;
note that the AS and IHX relations are closed in the subspace.

\begin{figure}[htpb]
$$ \psdraw{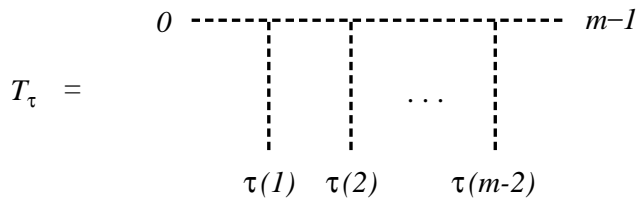}{3.6in} $$
\caption{The definition of $T_{\tau}$}\label{fig.Ttau}
\end{figure}

For an element $\tau$ in the symmetric group $\fS_{m-2}$
acting on the set $\{1,2,\cdots,m-2\}$,
let $T_{\tau} \in \cA(m)$ be the graph shown in Figure \ref{fig.Ttau}.

\begin{lem}\label{lem.basis}
We can take the set of $T_{\tau}$ as basis of the space $\cAt{m}$;
in particular the space is $(m-2)!$ dimensional.
\end{lem}

\begin{pf}
Let $D$ be any chord diagram in $\cAt{m}$.
We put red color on the path
connecting the two univalent vertices $0$ and $m-1$.
We can deform $D$ into a linear sum of $T_{\tau}$'s
by induction on the number of trivalent vertices
on the red path as follows.
We choose a trivalent vertex next to the red path,
and apply the IHX relation
regarding the segment connecting
the vertex and the red path as the character \lq\lq I''
in \lq\lq IHX''.
Then the number of trivalent vertices on the red path increases.
Hence we can show that
the space $\cAt{m}$ is spanned by the set of $T_{\tau}$.

In order to complete the proof of this lemma,
it is sufficient to prove that
$T_{\tau}$'s are linearly independent.
Suppose that $T_{\tau}$ could be expressed as
a linear sum of other $T_\s$'s.
We can \lq\lq substitute'' a Lie algebra $sl(m,\C)$
to dashed lines (see \cite{BarNatan})
to make a linear map of $\cAt{m}$ to $\C$.
For $j<k$ let $E_{jk}$ be the element in $sl(m,\C)$
which has $(j,k)$ entry $1$ and the other entries $0$.
We have a relation
$[E_{ij},E_{jk}]=E_{ik}$.
If we substitute $E_{12}$ to the univalent vertex $0$ and
$E_{k+1,k+2}$ to the vertex $\tau(k)$ for $k=1,2,\cdots,m-2$,
then $T_\s$ always vanishes unless $\s=\tau$,
though $T_{\tau}$ does not vanish
when we substitute the dual of $E_{1m}$ to the vertex $m-1$,
where we mean the dual with respect to the Killing form.
This is a contradiction, completing the proof.
\end{pf}

\subsection{Chord diagrams behaving
in a similar way as a solid circle}

We define $T_m \in \cAt{m}$ by
$$
T_m = \sum_{\tau\in{\fS}_{m-2}}
      \frac{(-1)^{r(\tau)}}{(m-1)\binom{m-2}{r(\tau)}} T_{\tau},
$$
where we denote by $r(\tau)$
the number of $k$ which satisfies $\tau(k)>\tau(k+1)$.
In this subsection,
our aim is to show Proposition \ref{prop.STUforTm}.
We begin with the following lemma,
which is a weak form of Proposition \ref{prop.STUforTm};
we can interchange any two adjacent univalent vertices
in Proposition \ref{prop.STUforTm}
whereas we can do it for particular pairs in Lemma \ref{lem.STUforT}.
We will show symmetries of $T_m$
in this subsection
to obtain Proposition \ref{prop.STUforTm} from Lemma \ref{lem.STUforT}.

\begin{figure}[htpb]
$$ \psdraw{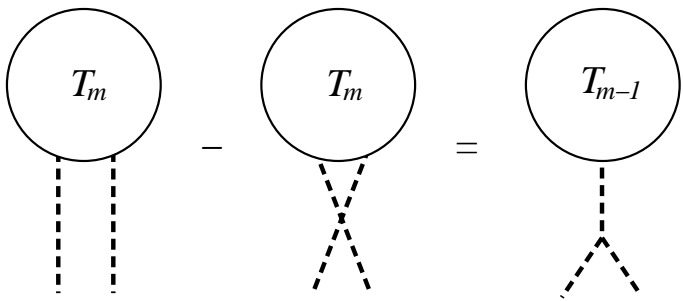}{3.6in} $$
\caption{A property of $T_m$ similar to the STU relation}\label{fig.STUinTm}
\end{figure}

\begin{lem}\label{lem.STUforT}
If $1 \le k \le m-3$,
then the difference between $T_m$ and
the chord diagram obtained from $T_m$
by changing two univalent vertices $k$ and $k+1$
can be expressed using $T_{m-1}$
as shown in Figure \ref{fig.STUinTm}.
\end{lem}

\begin{pf}
Since we can take the set of $T_{\tau}$ as basis of $\cAt{m}$,
we can express both sides of the required formula
as a linear sum of $T_{\tau}$'s.
It is sufficient to show that
the coefficients of $T_{\tau}$ in both sides
are equal for each $\tau$.

If $|\tau^{-1}(k)-\tau^{-1}(k+1)| \ge 2$,
then the coefficient of the left hand side
is equal to zero,
since $r(\tau)=r((k\ k+1)\circ \tau)$ holds in this case,
where we mean by $(k\ k+1)$ the interchange of $k$ and $k+1$.
On the other hand,
the coefficient of the right hand side is equal to zero,
since the right hand side is equal to a linear sum of $T_{\tau}$
for $\tau$ satisfying $|\tau^{-1}(k)-\tau^{-1}(k+1)|= \pm 1$;
we can see it by applying the IHX relation
in the right hand side.
Therefore the coefficients of $T_{\tau}$ in both sides
are equal in this case.

If $\tau^{-1}(k)-\tau^{-1}(k+1) = -1$,
then the coefficient of the left hand side is equal to
$t_{m,\tau}-t_{m,(k\ k+1)\circ \tau}$
where we put
$$
t_{m,\tau}= \frac{(-1)^{r(\tau)}}{(m-1)\binom{m-2}{r(\tau)}}.
$$
In this case we have
$r((k\ k+1)\circ \tau)=r(\tau)+1$
by the definition of $r(\cdot)$.
Hence we have
\begin{align*}
t_{m,\tau}-t_{m,(k\ k+1)\circ \tau}
&=\frac{(-1)^{r(\tau)}}{(m-1)\binom{m-2}{r(\tau)}}
 -\frac{(-1)^{r(\tau)+1}}{(m-1)\binom{m-2}{r(\tau)+1}} \\
&=\frac{(-1)^{r(\tau)}(m-r-2)}{(m-1)(m-2)\binom{m-3}{r(\tau)}}
 -\frac{(-1)^{r(\tau)+1}(r+1)}{(m-1)(m-2)\binom{m-3}{r(\tau)}}
 =\frac{(-1)^{r(\tau)}}{(m-2)\binom{m-3}{r(\tau)}}
\end{align*}
On the other hand,
the contribution of the right hand side to $T_{\tau}$
comes from $T_{\hat\tau} \in \cA(m-1)$,
where we define $\hat\tau \in \fS_{m-1}$ by
putting $\hat\tau^{-1}(j)$ to be
$\tau^{-1}(j)$ if $j \le k$, $\tau^{-1}(j+1)-1$ if $j > k$.
The coefficient of $T_{\hat\tau}$ is equal to $t_{m-1,\hat\tau}$;
in this case we have
$r(\hat\tau)=r(\tau)$
by the definition of $r(\cdot)$.
Therefore the coefficients of both sides are equal,
completing this case.

If $\tau^{-1}(k)-\tau^{-1}(k+1) = 1$,
we can show that the coefficients are equal
in a similar way as above,
completing the proof.
\end{pf}

\begin{figure}[htpb]
$$ \psdraw{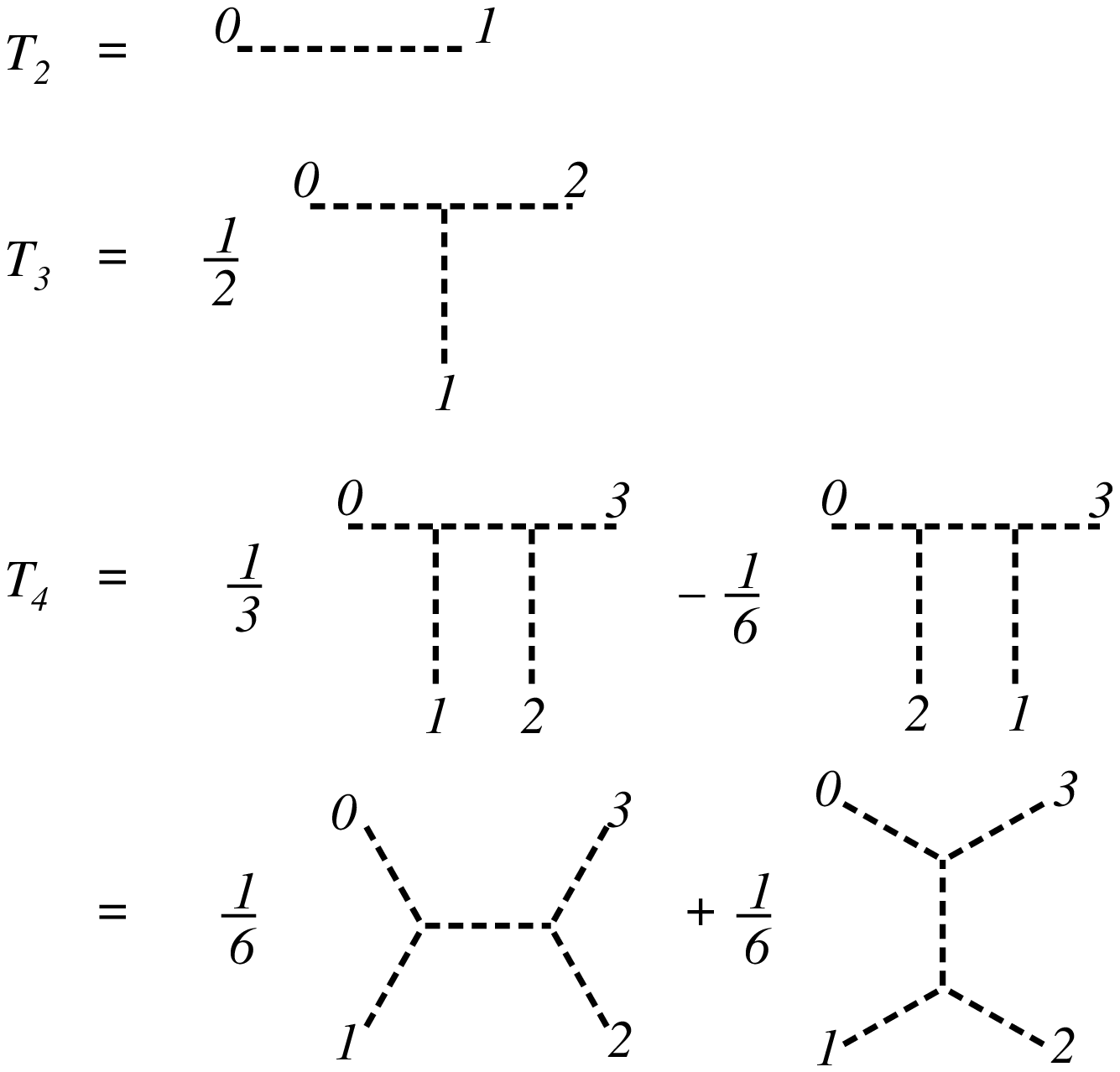}{3.6in} $$
\caption{Simple cases of $T_m$}\label{fig.T2T3T4}
\end{figure}

\begin{lem}\label{lem.sym1}
The chord diagram $T_m$ is symmetric
with respect to mirror image
which replaces
$0,1,\cdots,m-1$ with $m-1,m-2,\cdots,0$ respectively,
see Figure \ref{fig.T2T3T4} for simple cases.
To be exact, the inversion
replaces $T_m$ to $(-1)^{m-2}T_m$
because it changes the orientations of $m-2$ trivalent vertices.
\end{lem}

\begin{pf}
By the inversion which replaces
$0,1,\cdots,m-1$ with $m-1,m-2,\cdots,0$,
the chord diagram $T_{\tau}$ moves to $(-1)^{m-2}T_{\tau'}$ where
$$
\tau' =
\begin{pmatrix}
1 & 2 & \cdots & m-2 \\ m-2 & m-1 & \cdots & 1
\end{pmatrix}
\circ \tau \circ
\begin{pmatrix}
1 & 2 & \cdots & m-2 \\ m-2 & m-1 & \cdots & 1
\end{pmatrix}
$$
and the sign is derived from
the number of trivalent vertices;
we use the AS relation such times.
We can obtain $r(\tau)=r(\tau')$
by definition of $r(\cdot)$.
Therefore the inversion maps $T_m$ to $(-1)^{m-2} T_m$,
completing the proof.
\end{pf}

\begin{lem}\label{lem.sym2}
The chord diagram $T_m$ is symmetric
with respect to mirror image
which replaces
$0,1,2,\cdots,m-2,m-1$ with $m-2,m-3,m-4,\cdots,0,m-1$ respectively.
To be exact,
the inversion replaces $T_m$
to $(-1)^{m-2}T_m$
because of the orientations of trivalent vertices.
\end{lem}

\begin{pf}
Consider the linear map $i$
of $\cAt{m}$ to $\cAt{m+1}$
which maps a chord diagram $D$
to the diagram $D$ added a small branch
near the univalent vertex $m-1$
as shown in Figure \ref{fig.addvertex}.
The map $i$ is an injection;
we can see it
by taking basis of $\cAt{m}$ and $\cAt{m+1}$
in the way how we choose the red path
connecting $0$ and $m-1$ and
connecting $0$ and $m$ respectively
as in the proof of Lemma \ref{lem.basis}.

\begin{figure}[htpb]
$$ \psdraw{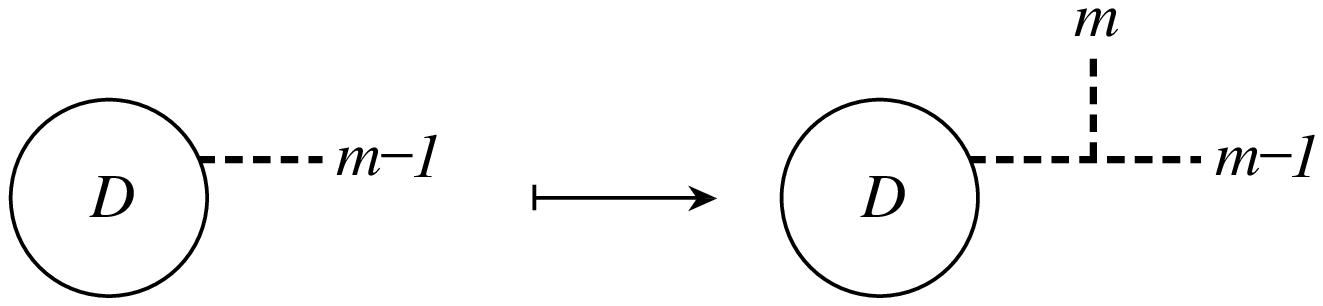}{3.6in} $$
\caption{The map of $\cAt{m}$ to $\cAt{m+1}$}\label{fig.addvertex}
\end{figure}

We denote by $\fS'_{m-1}$
the symmetric group acting on the set
$\{0,1,\cdots,m-2\}$.
For an element $\s \in \fS'_{m-1}$,
let $S_\s \in \cAt{m+1}$ be the chord diagram
shown in Figure \ref{fig.Ssigma}.
By Lemma \ref{lem.basis}
we can take the set of $S_\s$ as basis of $\cAt{m+1}$.

\begin{figure}[htpb]
$$ \psdraw{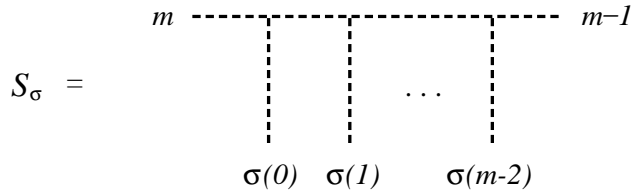}{3.6in} $$
\caption{The definition of $S_\s$}\label{fig.Ssigma}
\end{figure}

We can express $i(T_m)$ as in Lemma \ref{lem.TtoS} below.
In order to complete the proof of the present lemma,
it is sufficient to show that
$i(T_m)$ is symmetric
with respect to the inversion
replacing $0,1,\cdots,m-2$ with $m-2,m-1,\cdots,0$,
since the map $i$ is an injection.
The inversion maps $S_\s$ to $S_{\s'}$ where
$$
\s' =
\begin{pmatrix}
0 & 1 & \cdots & m-2 \\ m-2 & m-1 & \cdots & 0
\end{pmatrix}
\circ \s,
$$
note that we have no change of sign in this case
because we fix the ends $m-1$ and $m$ of the \lq\lq red path''.
We can obtain
$r(\s')=m-2-r(\s)$
by the definition of $r(\cdot)$.
By Lemma \ref{lem.TtoS} below,
the inversion maps $i(T_m)$ to $(-1)^{m-2} i(T_m)$,
completing the proof.
\end{pf}

\begin{lem}\label{lem.TtoS}
$$
i(T_m) =
    \sum_{\s\in{\fS'}_{m-1}}
      \frac{(-1)^{r(\s)}}{(m-1)\binom{m-2}{r(\s)}} S_\s.
$$
\end{lem}

\begin{pf}
We put $S_{m+1}=i(T_m)$.
Since the set of $S_\s$ is basis of $\cAt{m+1}$,
we can put $S_{m+1}=\sum_\s s_\s S_\s$
with some scalars $s_\s$.
In the following of this proof
we will show the following formula
\begin{equation}\label{eq.snodef}
s_\s = \frac{(-1)^{r(\s)}}{(m-1)\binom{m-2}{r(\s)}}
\end{equation}
in two steps by induction on $m$.

\noindent
{\bf Step 1.}\quad
If $\s$ is a cyclic permutation
$(0,1,2,\cdots,k)$ for an integer $k$,
we can obtain (\ref{eq.snodef}) by definition of $i(T_m)$ as follows.
We consider that
what $T_{\tau}$ contributes to such $S_\s$
after changing basis of $\cAt{m}$ to that of $\cAt{m+1}$;
we expand the left picture in Figure \ref{fig.expandT}
using the IHX relation to obtain the right picture.
We use the IHX relation
replacing \lq\lq I'' with the difference of \lq\lq H'' and \lq\lq X''.
Noting that the univalent vertex $0$
interchanges $k$ times with another vertex,
we must use \lq\lq X'' $k$ times in the expansion.
Hence the possibilities of $\tau$ are as follow;
\begin{align*}
\tau &=
\left( {\begin{matrix}
1 &2 &\cdots &l_k-1 &l_k &l_k+1 &\cdots \\
k+1 &k+2 &\cdots &l_k+k-1 &k &l_k+k &\cdots
\end{matrix}\quad
\begin{matrix}
l_1-1 &l_1 &l_1+1 &\cdots &m-2 \\
l_1 &1 &l_1+1 &\cdots &m-2
\end{matrix}} \right) \\
&\text{or }
\left( {\begin{matrix}
1 &2 &\cdots &l_{k-1}-1 &l_{k-1} &l_{k-1}+1 &\cdots \\
k &k+1 &\cdots &l_{k-1}+k-2 &k-1 &l_{k-1}+k-1 &\cdots
\end{matrix}\quad
\begin{matrix}
m-2 \\
m-2
\end{matrix}} \right)
\end{align*}
In the first type
we can freely choose $\{ l_k,l_{k-1},\cdots,l_1 \}$
from $\{2,3,\cdots,m-2 \}$.
Hence there are $\binom{m-3}{k}$ possibilities of $\tau$,
and $r(\tau)=k$ holds in this type.
Similarly we have $\binom{m-3}{k-1}$ possibilities of $\tau$,
which satisfy $r(\tau)=k-1$ in the second type.
Therefore we have
$$
s_\s = (-1)^k \binom{m-3}{k}
       \frac{(-1)^k}{(m-1)\binom{m-2}{k}}
     + (-1)^k \binom{m-3}{k-1}
       \frac{(-1)^{k-1}}{(m-1)\binom{m-2}{k-1}}
$$
where the term $(-1)^k$ is derived from
the number of usage of \lq\lq X''.
This formula satisfies (\ref{eq.snodef}), completing Step 1.

\begin{figure}[htpb]
$$ \psdraw{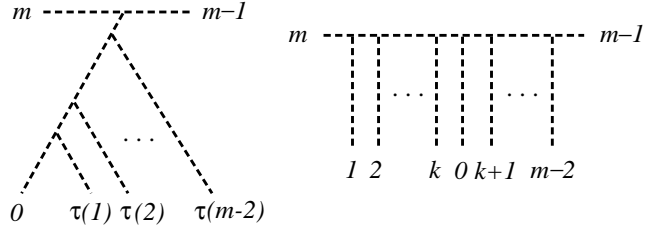}{3.6in} $$
\caption{Expanding $T_{\tau}$}\label{fig.expandT}
\end{figure}

\noindent
{\bf Step 2.}\quad
In this step we will show that
if (\ref{eq.snodef}) holds for $\s$
then it also holds for $(k\ k+1)\circ \s$
for any $k = 1,2,\cdots,m-3$,
where we mean by $(k\ k+1)$ the interchange of $k$ and $k+1$.

We have the formula shown in Figure \ref{fig.STUforS},
where the first and third equalities in the figure
are derived from the definition of $S_\star$
and the second equality is derived from Lemma \ref{lem.STUforT}.

\begin{figure}[htpb]
$$ \psdraw{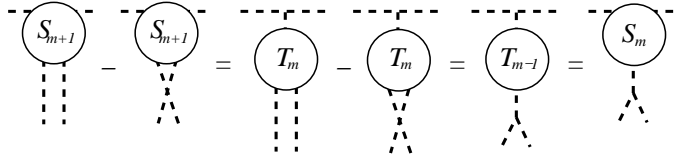}{3.6in} $$
\caption{$S_\star$ satisfies a relation similar to the STU relation}
\label{fig.STUforS}
\end{figure}

Hence $S_{m+1}$ satisfies a relation similar to the STU relation;
this means
$$
s_\s - s_{(k\ k+1) \circ \s} =
\begin{cases}
s_{\hat\s} \quad & \text{ if } \s^{-1}(k)-\s^{-1}(k+1)=-1 \\
-s_{\hat\s} \quad & \text{ if } \s^{-1}(k)-\s^{-1}(k+1)=1 \\
0 \quad & \text{ if } |\s^{-1}(k)-\s^{-1}(k+1)| \ge 2
\end{cases}
$$
where we define $\hat\s$ by putting
$\hat\s^{-1}(j)$ to be $\s^{-1}(j)$ if $j \le k$,
$\s^{-1}(j+1)-1$ of $j>k$;
note that
we can obtain the chord diagram $S_{\hat\s}$
from $S_\s$ by gluing two adjacent dashed edges
who have univalent vertices $k$ and $k+1$ respectively,
and we can define $\hat\s$
only when $\s^{-1}(k)-\s^{-1}(k+1)=\pm1$.
We note that we can use (\ref{eq.snodef}) for $\hat\s$
by hypothesis of induction.
Since we can check that the above formula satisfies (\ref{eq.snodef}),
we obtain the required claim of Step 2.

Since we can obtain any $\s \in \fS_{m-1}$
from some permutation $(0,1,2,\cdots,k)$
by composing interchanges $(j\ j+1)$,
we can show (\ref{eq.snodef}) for any $\s$,
completing the proof.
\end{pf}

By Lemmas \ref{lem.sym1} and \ref{lem.sym2},
we immediately obtain the following proposition.

\begin{prop}\label{prop.sym}
The chord diagram $T_m$ is symmetric
with respect to the action of dihedral group of order $2m$;
we show simple cases in Figure \ref{fig.T2T3T4}.
Note that the sign of $T_m$ possibly changes by the action
as in the statements of Lemmas \ref{lem.sym1} and \ref{lem.sym2}.
\end{prop}

Noting the symmetry of $T_m$ in Proposition \ref{prop.sym},
we can obtain the following proposition
from Lemma \ref{lem.STUforT}.

\begin{prop}\label{prop.STUforTm}
The difference between $T_m$ and the chord diagram obtained by changing
any adjacent two univalent vertices is equal to
$T_{m-1}$ with one extra trivalent vertex
as shown in Figure \ref{fig.STUinTm}.
\end{prop}

\subsection{Replacing solid circles with dashed lines}

\begin{figure}[htpb]
$$ \psdraw{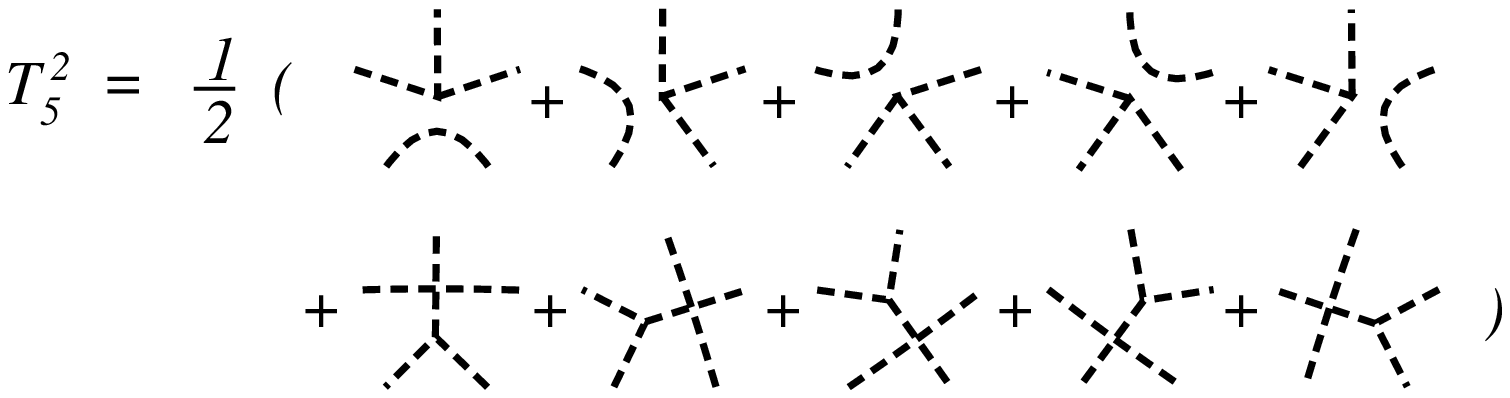}{3.6in} $$
\caption{The picture of $T^2_5$}\label{fig.T25}
\end{figure}

Let $n$ be a positive integer.
For an integer $m$ with $m \ge 2n$,
we define $T^n_m \in \cA(m)$ as follows.
We divide $m$ into $n$ integers as;
$$
m=m_1+m_2+\cdots+m_n, \qquad m_1 \ge m_2 \ge \cdots \ge m_n \ge 2.
$$
Let $T^n_m$ be the sum of
all configurations of disjoint union of $T_{m_i}$'s
where we take the configurations preserving
cyclic order of univalent vertices of each $T_{m_i}$.
We show the picture of $T^2_5$ in Figure \ref{fig.T25}.
If $m$ is less than $2n$,
then we put $T^n_m$ to be $0$.

\begin{prop}\label{prop.STUforTnm}
For any positive integer $n$,
the difference between $T^n_m$ and the chord diagram obtained by changing
any adjacent two univalent vertices is equal to
$T^n_{m-1}$ with one extra trivalent vertex.
Namely $T^n_m$ also satisfies the same formula
as in Figure \ref{fig.STUinTm} for $T_m$.
\end{prop}

\begin{pf}
For each term in $T^n_m$,
we consider the difference in the left hand side.
If the two adjacent vertices are
on different connected components of $T^n_m$,
then the difference vanishes.
If they are on the same connected component,
then the difference is equal to a term in $T^n_{m-1}$
in the right hand side by Proposition \ref{prop.STUforTm}.
Hence we obtain the required formula.
\end{pf}

\begin{figure}[htpb]
$$ \psdraw{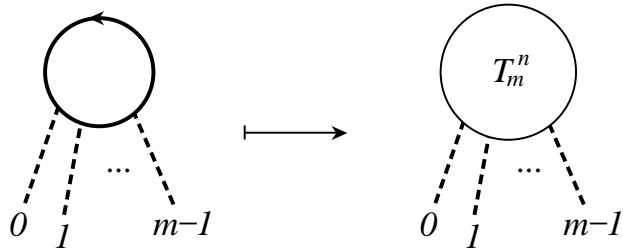}{3.6in} $$
\caption{The definition of $\iota_n$}\label{fig.iota}
\end{figure}

Using the above proposition,
we can define a linear map
$$
\iota_n:\cA(\coprod^l S^1)
\longrightarrow\cA(\phi)
$$
by putting
$\iota_n(\text{a solid circle with $m$ dashed univalent vertices})$
to be $T^n_m$ as shown in
Figure \ref{fig.iota}.
We also denote by the same notation $\iota_n$
the naturally extended map of $\hcA(\coprod^l S^1)$
to $\hcA(\phi)$.

\section{A series of invariants of 3-manifolds}

In this section
we will construct a series of
topological invariants $\Omega_n(M)$
of a 3-manifold $M$.
We show the invariance of $\cZ(L)$
under orientation change and Kirby move II
using the equivalence relation $P_{\star}$ defined below.
Since the relation $P_{n+1}$ vanish
in low degrees in the image of $\iota_n$,
we will obtain a series of invariants
in low degrees of $\cA(\phi)$ through the maps $\iota_n$.

\subsection{Invariance under orientation change and Kirby move II}

\begin{figure}[htpb]
$$ \psdraw{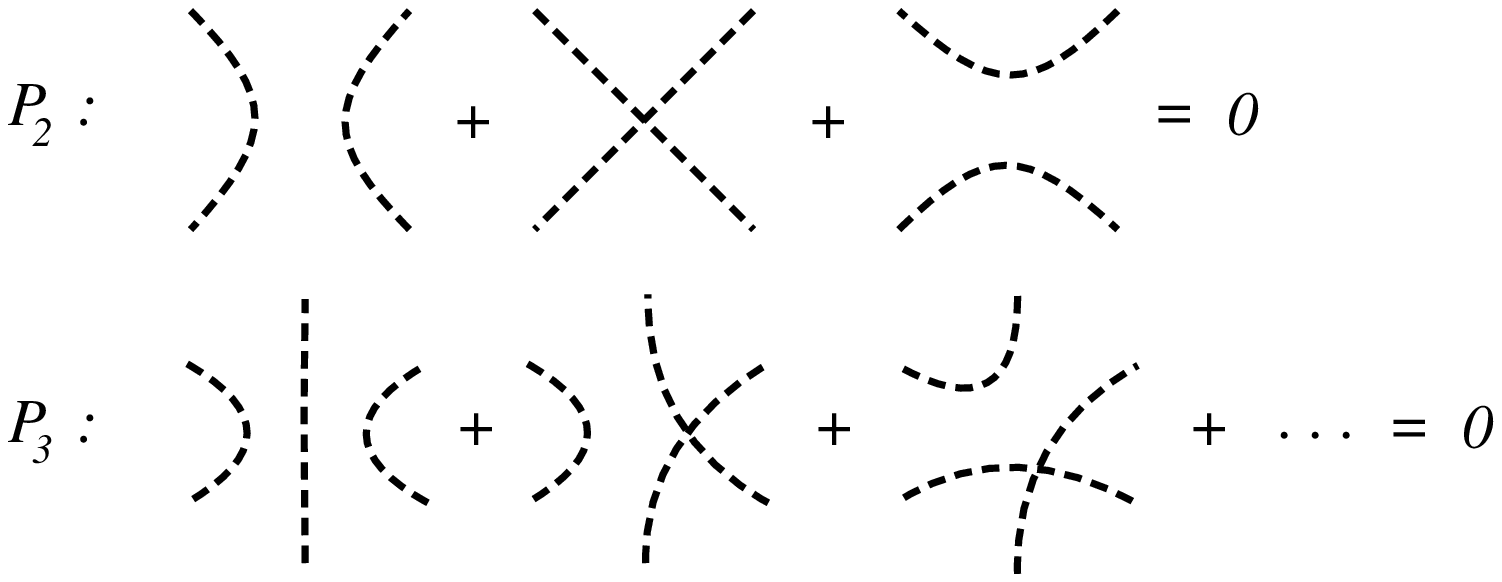}{3.6in} $$
\caption{The relations $P_2$ and $P_3$}\label{fig.P2P3}
\end{figure}

Let $\oA(X)$ be the space of chord diagrams including
dashed trivial circles.
We denote by $\hoA{X}$ its completion.
For each positive integers $n$,
we define an equivalence relation $P_n$
in $\oA(X)$ and $\hoA{X}$ as follows.
Let $P_1$ be the equivalence relation such that
any chord diagram with non-empty dashed lines is equivalent to zero.
Let $P_2$ be the equivalence relation shown in Figure \ref{fig.P2P3};
the left hand side of the first formula in Figure \ref{fig.P2P3}
is the sum over all pairings of $4$ points.
Similarly we define the equivalence relation $P_n$
such that the sum over all pairings of $2n$ points is equivalent to zero.

Before proving the invariance of $\cZ(L)$
in Proposition \ref{prop.invariantbyP},
we prepare the following lemma.

\begin{lem}\label{lem.decreaselegs}
Let $D$ be a chord diagram on $X$,
and $C$ a component of $X$.
\begin{enumerate}
\item
Then, with the equivalence relation $P_{n+1}$,
the chord diagram $D$ is equivalent to
a linear sum of chord diagrams
each of which has at most $2n$ univalent vertices on $C$.
\item
Further, the chord diagram $D$ becomes equivalent to
a linear sum of chord diagrams
each of which has either $n$ isolated dashed chords on $C$
or at most $2n-1$ univalent vertices on $C$.
Here we mean by an isolated chord
a dashed arc with no trivalent vertices and
two adjacent univalent vertices on $C$.
\end{enumerate}
\end{lem}

\begin{figure}[htpb]
$$ \psdraw{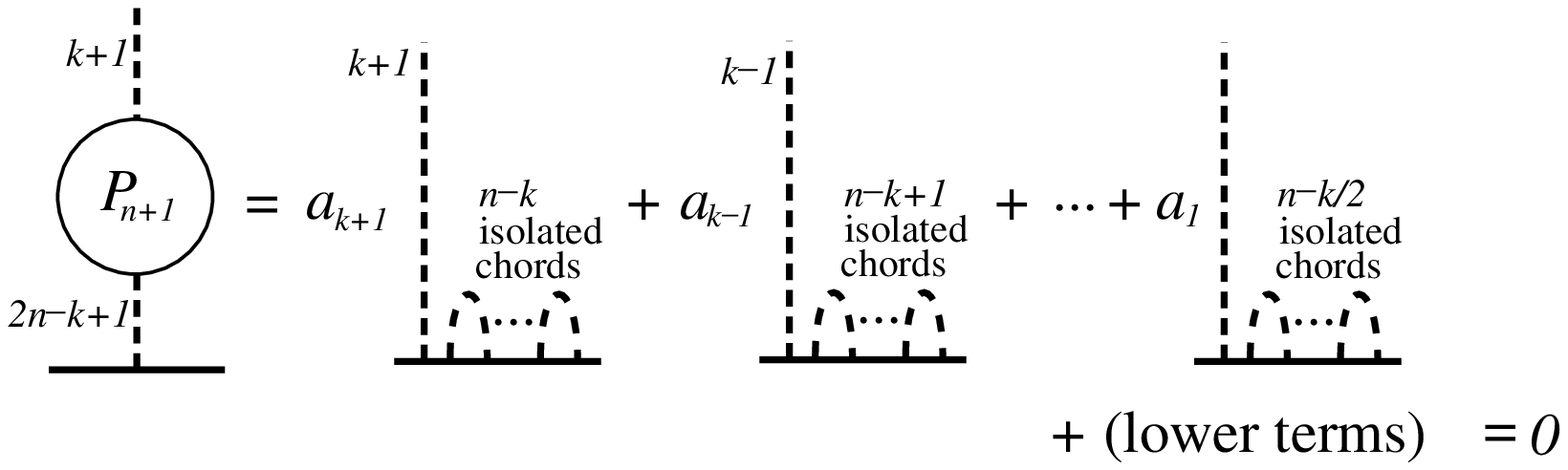}{3.6in} $$
\caption{Using $P_{n+1}$ to decrease univalent vertices;
this picture is for even $k$}\label{fig.useP}
\end{figure}

\begin{pf}
We will see the case $n=3$ before the general case.
It is sufficient to show that,
if the diagram $D$ has $k$ univalent vertices on $C$ with $k>4$,
then it is equivalent to a linear sum of
chord diagrams each of which has
less than $k$ univalent vertices;
we will call them lower terms.
We use the relation $P_3$ as in Figure \ref{fig.useP3},
where the second equality is derived from the STU relation;
we can use the STU relation to interchange
two univalent vertices modulo a lower term.

\begin{figure}[htpb]
$$ \psdraw{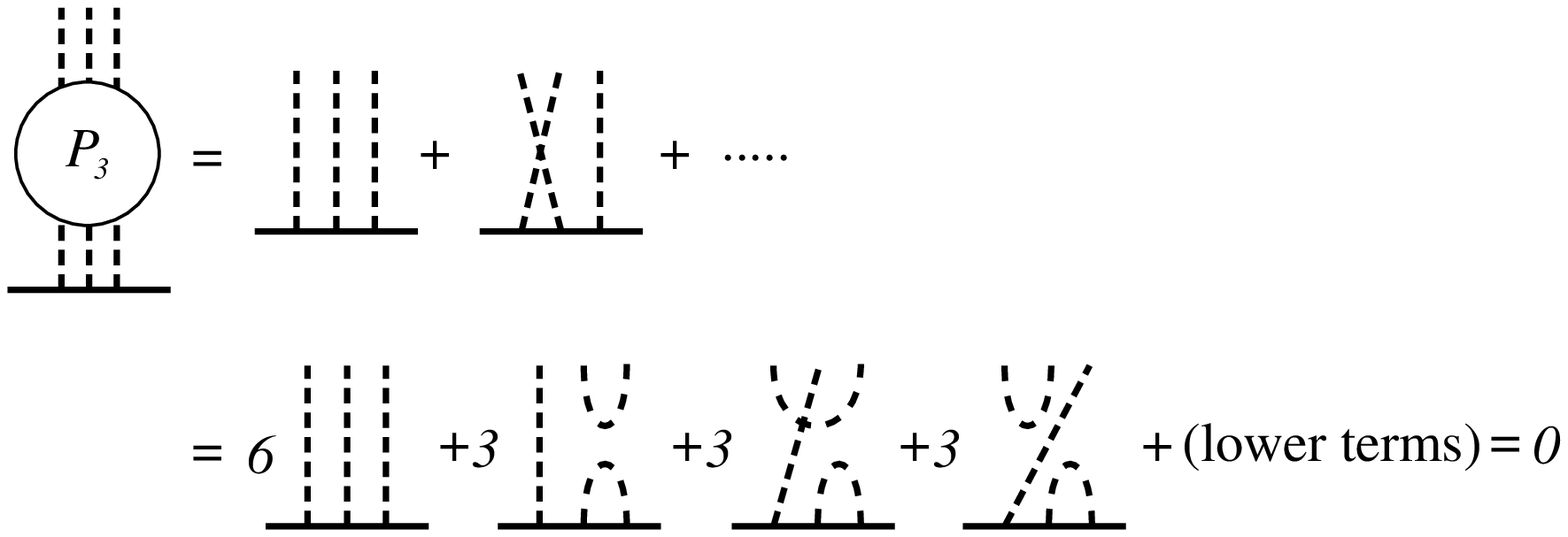}{3.6in} $$
\caption{Using $P_3$ to make isolated chords}\label{fig.useP3}
\end{figure}

Then we can replace $D$ with
a chord diagram with an isolated chord,
where we mean by an isolated chord a trivial dashed arc
with adjacent univalent vertices on $X$.
Iterating this procedure,
we can replace $D$ with
a chord diagram with at least two isolated chords.
Then we can use the relation $P_3$ again
as in Figure \ref{fig.useP3again},
and we can replace the diagram with lower terms,
completing the proof of (1) for the case $n=3$.

\begin{figure}[htpb]
$$ \psdraw{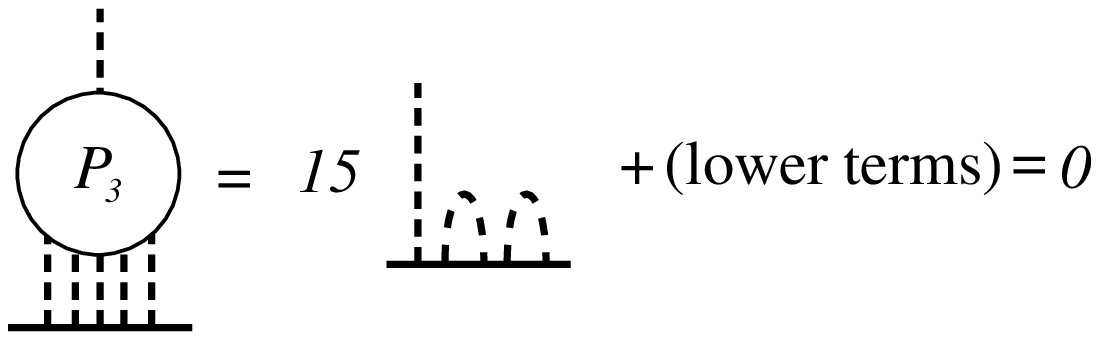}{3.6in} $$
\caption{Using $P_3$ to decrease univalent vertices}\label{fig.useP3again}
\end{figure}

In order to prove (2),
it is sufficient to show that
any chord diagram with $4$ univalent vertices on $C$
is equivalent to a chord diagram with
two isolated chords on $C$ modulo lower terms.
We use the relation in Figure \ref{fig.useP3} again,
to make one isolated chord.
We further use the relation $P_3$
as in Figure \ref{fig.useP3further}.
Then we can replace
the diagram with a diagram with two isolated chords on $C$,
completing the proof of (2) for the case $n=3$.

\begin{figure}[htpb]
$$ \psdraw{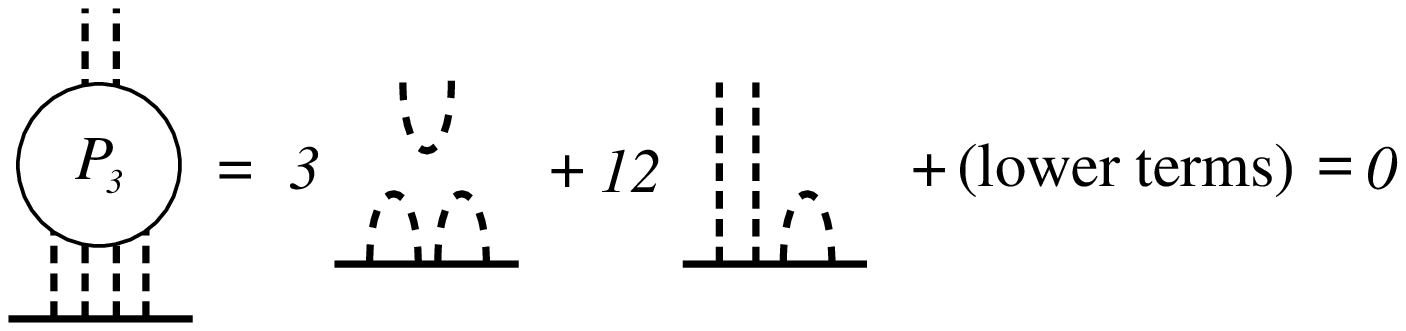}{3.6in} $$
\caption{Using $P_3$ to make two isolated chords}\label{fig.useP3further}
\end{figure}

In a general case for proving (1),
it is sufficient to show that,
if the diagram $D$ has $m$ univalent vertices on $C$ with $m>2n$,
then it is equivalent to
a linear sum of chord diagrams
with at most $m-1$ univalent vertices.
We use the relation $P_{n+1}$ as in Figure \ref{fig.useP}
for $0 \le k \le n$;
we use it for $k=n$ at the beginning to make an isolated chord,
and use it for $k=2,3,\cdots$ to increase isolated chords,
and finally use it for $k=0$
to replace the diagram with lower terms;
note that we can do it assuming $m>2n$.
This completes the proof of (1).

In order to prove (2),
it is sufficient to show that
any chord diagram with $2n$ univalent vertices on $C$
is equivalent to a diagram
with $n$ isolated chords on $C$ modulo lower terms.
We use the relation in Figure \ref{fig.useP} as above
to make $n-1$ isolated chords on $C$.
We further use the relation for $k=2$,
and we obtain a chord diagram with $n$ isolated chords on $C$.
This completes the proof of (2).
\end{pf}

\begin{rem}
We can uniquely characterize $T_m^n$
as an element of $\cA(m)/P_{n+1}$
by the following four conditions.
\begin{enumerate}
\item
It satisfies $T_m^n=0$ if $m<2n$.
\item
It is invariant under cyclic permutation
of the $m$ external univalent vertices.
\item
It satisfies a relation similar to the STU relation;
in other words, satisfies Proposition \ref{prop.STUforTnm}.
\item
It contains neither dashed cycle nor dashed component; that is,
it is equal to a linear sum of chord diagrams
each of which is a disjoint union of
simply connected chord diagrams having external vertices.
\end{enumerate}
Outline of a proof of the fact is as follows.
By Lemma \ref{lem.decreaselegs}
a solid circle with $m$ dashed chords
is equivalent
(modulo $P_{n+1}$ and lower terms)
to the disjoint union of
a solid circle with $n$ isolated dashed chords
and a linear sum of dashed trivalent graphs.
We can obtain an alternative definition
of $\iota_n$;
we define the inverse map $\iota_n^{-1}$
by removing the solid circle with $n$ dashed isolated chords
from the above disjoint union.
The image of solid circle with $m$ short dashed chords
through $\iota^{-1}_n$
is equal to $T_m^n$ by definition of $\iota_n$, and
we can check the image satisfies the above four conditions.
\end{rem}

By Lemma \ref{lem.decreaselegs}
we can replace a solid circle
with a solid circle with $2n$ univalent vertices
modulo lower terms.
In the following proposition
we show the invariance of $\cZ(L)$
by reducing the proof to the above particular solid circle
ignoring the lower terms.

\begin{prop}\label{prop.invariantbyP}
Let $L$ be any oriented framed link of $l$ components,
$n$ any positive integer.
\begin{enumerate}
\item
The equivalence class $[\cZ(L)]$ including $\cZ(L)$
in $\hoA{\coprod^l S^1}/L_{<2n}, P_{n+1}, O_n$
does not depend on orientation of $L$,
where we denote by $L_{<2n}$
the equivalence relation such that
any chord diagram including a solid circle
with less than $2n$ dashed univalent vertices is equivalent to zero,
and $O_n$ the equivalence relation such that
a trivial dashed circle is equivalent to $-2n$.
\item
The above equivalence class
is invariant under Kirby move II.
\end{enumerate}
\end{prop}

\begin{pf}
We will show the proposition for each chord diagram $D$ composing $\cZ(L)$.

Let $L'$ be the framed link obtained from $L$ by
changing the orientation of a component $C$;
we denote by $C'$
the corresponding solid circle in the chord diagram in $\cZ(L)$.
The change of $\cZ(L)$ to $\cZ(L')$
is given in (\ref{eq.revori});
the diagram in $\cZ(L')$ obtained by changing $D$
is equal to $S_{(C)}(D)$.
We will show that $S_{(C)}(D)$ is equivalent to $D$.

If $D$ has less than $2n$ dashed univalent vertices on $C'$,
then both of $D$ and $S_{(C)}(D)$ vanish
by the relation $L_{<2n}$.
Hence they are equivalent.

If $D$ has $2n$ univalent vertices on $C'$,
then by (\ref{eq.revori})
we can obtain $S_{(C)}(D)$ from $D$
by changing the order of univalent vertices on $C'$.
By the relation $L_{<2n}$ and the STU relation,
their equivalence classes are invariant
under changing the order of univalent vertices on $C'$
in this case.
Hence the equivalence classes are equivalent.

If $D$ has $k$ univalent vertices on $C'$ for $k>2n$,
then we can reduce this case to the case $k=2n$ as follows.
By Lemma \ref{lem.decreaselegs}
we can replace $D$ with a chord diagram with less univalent vertices on $C$.
The relations we used in the proof are
only the STU relation and the relation $P_{n+1}$.
Namely, we have a series of a linear sum of chord diagrams;
$D=D_0,D_1,D_2,\cdots D_N$,
where $D_i$ and $D_{i+1}$ are related by
either of the STU relation or the relation $P_{n+1}$.
Note that, if $D_i$ and $D_{i+1}$ are related by the STU relation,
then $S_{(C)}(D_i)$ and $S_{(C)}(D_{i+1})$ are also related by the relation,
because $S_{(C)}$ and the STU relation commute
as shown in Figure \ref{fig.SandSTU};
in fact, it guarantees that
the map $S_{(C)}$ is well defined.
Further note that
the relation $P_{n+1}$ and the STU relation commute,
because the relation $P_{n+1}$ is independent of solid chords.
Hence $S_{(C)}(D_i)$ and $S_{(C)}(D_{i+1})$ are
related by the STU relation or the relation $P_{n+1}$.
Therefore, if $D_{i+1}$ is equivalent to $S_{(C)}(D_{i+1})$,
then if $D_i$ is equivalent to $S_{(C)}(D_i)$;
this implies the required reduction,
completing this case.
This completes the proof of (1).

\begin{figure}[htpb]
$$ \psdraw{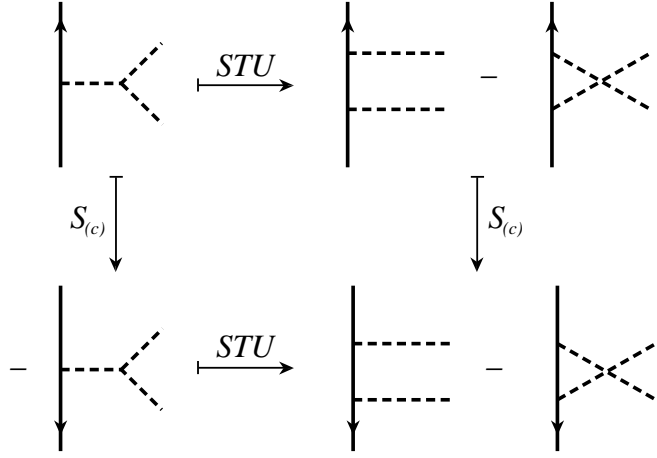}{3.6in} $$
\caption{The map $S_{(C)}$ and the STU relation commute}\label{fig.SandSTU}
\end{figure}

Let $L''$ be a framed link obtained from $L$
by Kirby move II
taking handle sliding of some component over a component $C$.
The change of $\cZ(L)$ to $\cZ(L'')$
is given by Proposition \ref{prop.KII};
the diagram in $\cZ(L'')$ obtained by changing $D$
is equal to $\D_{(C)}(D)$.
We will show that $\D_{(C)}(D)$ is equivalent to $D$.

If there are at most $2n-1$ univalent vertices on $C'$,
then both of $D$ and $\D_{(C)}(D)$ vanish
by the relation $L_{<2n}$.
Hence they are equivalent.

If there are $2n$ univalent vertices on $C'$,
then by Proposition \ref{prop.KII}
$\D_{(C)}(D)$ is equal to a sum of
$D$ and the other terms.
Further there are less than $2n$ univalent vertices on $C'$
for each of the other terms.
Hence their equivalent classes vanish by the relation $L_{<2n}$.
Therefore $D$ and $\D_{(C)}(D)$ are equivalent.

If there are $k$ univalent vertices on $C'$ for $k>2n$,
then we can reduce this case to the case $k-1$
as in the above proof of (1).
Instead of the commutation of $S_{(C)}$ and the STU relation,
we need the commutation of $\D_{(C)}$ and the STU relation
in the present case; we show it in Figure \ref{fig.DandSTU}.

\begin{figure}[htpb]
$$ \psdraw{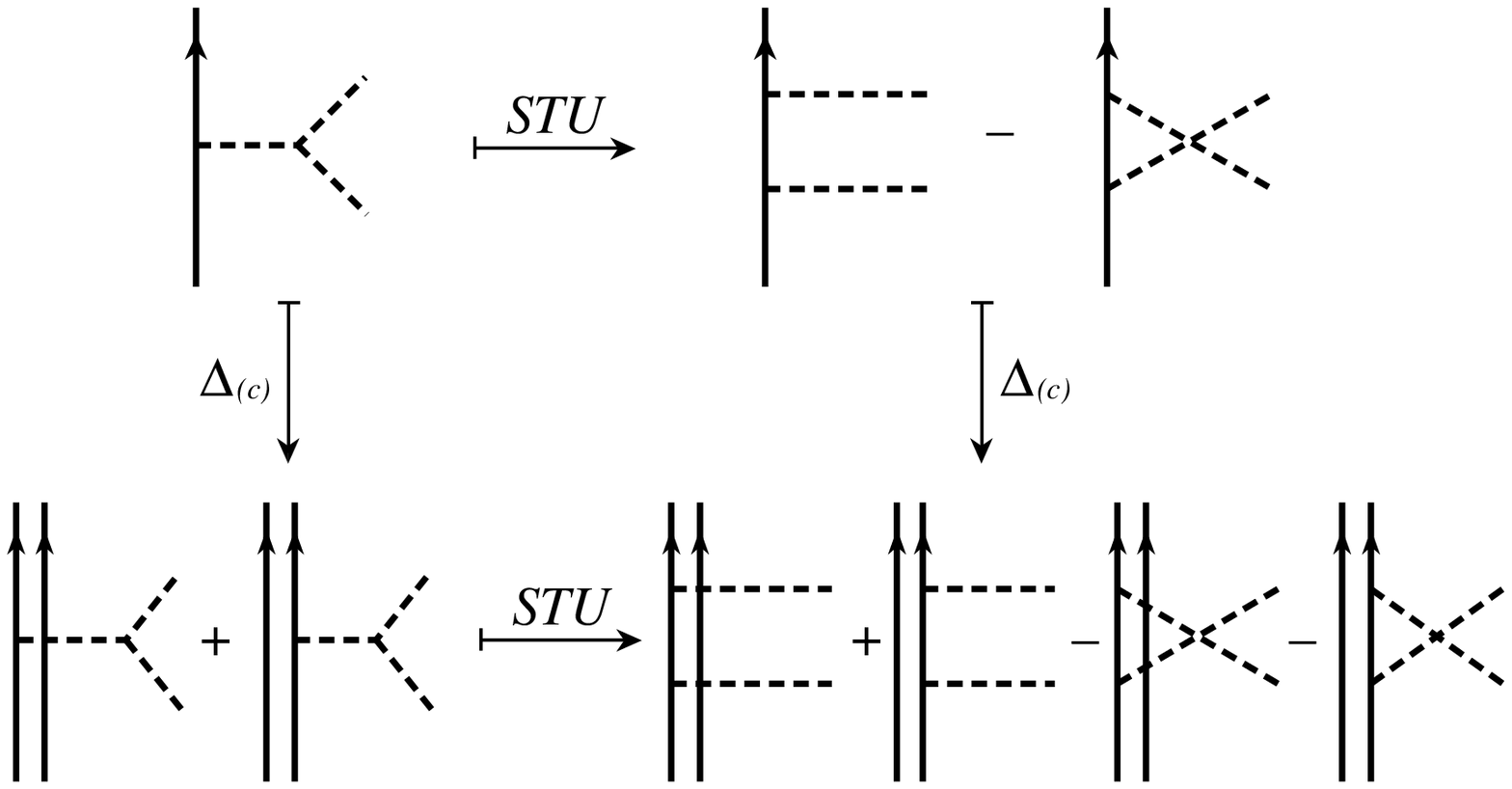}{3.6in} $$
\caption{The map $\D_{(C)}$ and the STU relation commute}\label{fig.DandSTU}
\end{figure}

This completes the proof of Proposition \ref{prop.invariantbyP}.
\end{pf}

\subsection{Moving $\cZ(L)$ into a set with algebra structure}

In order to show the invariance under Kirby move I,
we move $\cZ(L)$ into a quotient space of $\cA(\phi)$,
in which there is an algebra structure
with respect to the disjoint union of chord diagrams.
Before moving, we prepare the following lemma,
which guarantees that
we need not consider $P_{n+1}$
in low degrees of $\cA(\phi)$.

\begin{lem}\label{lem.Pvanishes}
The identity map induces an isomorphism of
the quotient space $\cA(\phi)/D_{>n}$
to the quotient space $\oA(\phi)/D_{>n},P_{n+1},O_n$.
\end{lem}

\begin{pf}
We can remove the trivial dashed component in $\oA(\phi)$ by $O_n$.
Hence it is sufficient to show the claim that,
if an element of $\oA(\phi)$ including $P_{n+1}$,
then either it vanishes or its degree is greater than $n$,
where we also denote by $P_{n+1}$
the formula defining the relation $P_{n+1}$.
Note that $P_{n+1}$ is symmetric
with respect to the action of $\fS_{2n}$
acting on the set of $2n$ ends of $P_{n+1}$.
We will show the claim for any outside of $P_{n+1}$.

We will show that the degree must be more than $n$
unless the element vanishes.
We can make an injection of the set of $2n+2$ ends of $P_{n+1}$
to the set of trivalent vertices in the outside as follows.
Choose an end of $P_{n+1}$.

If the end is not connected to any other end in the outside,
we associate one of trivalent vertices
in the connected component of the end in the outside.
Otherwise we can find a path
connecting the end to one of the other ends in the outside.

If there are no trivalent vertices on the path,
the element vanishes with the relation $O_n$
as shown in Figure \ref{fig.Pvanish}.

\begin{figure}[htpb]
$$ \psdraw{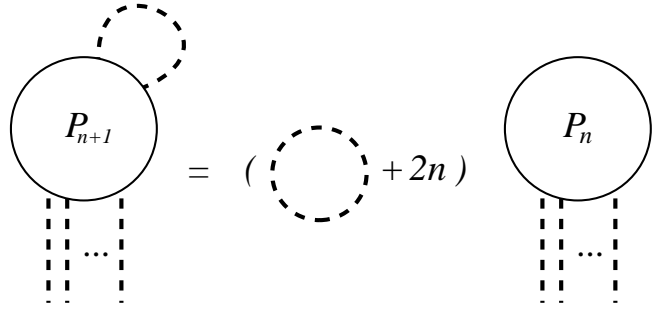}{3.6in} $$
\caption{The relation $P_{n+1}$ vanishes}\label{fig.Pvanish}
\end{figure}

If there are one trivalent vertices on the path,
the element vanishes using the symmetry of $P_{n+1}$
and the AS relation
as shown in Figure \ref{fig.Pvanishagain}.

\begin{figure}[htpb]
$$ \psdraw{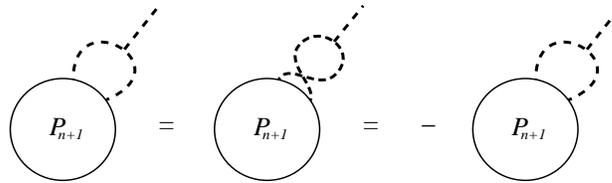}{3.6in} $$
\caption{The relation $P_{n+1}$ vanishes again}\label{fig.Pvanishagain}
\end{figure}

Otherwise there must be at least two trivalent vertices on the path.
Then we associate the nearest trivalent vertex to the end,
to obtain an injection of the set of $2n+2$ ends of $P_{n+1}$
to the set of trivalent vertices.

Hence we showed that
the number of trivalent vertices is greater than $2n$;
this implies that the degree is greater than $n$,
completing the proof.
\end{pf}

\begin{prop}\label{prop.invariant}
\begin{enumerate}
\item
The equivalence class
$[\iota_n(\cZ(L))] \in \cA(\phi)/D_{>n}$
does not depend on orientation of $L$.
\item
Further, the above equivalence class
does not change under Kirby move II.
\end{enumerate}
\end{prop}

\begin{pf}
By the map $\iota_n$
any chord diagram with less than $2n$ univalent vertices
on a solid circle of the diagram vanishes.
Hence we have the following map;
$$
\hoA{\coprod^l S^1}/L_{<2n},P_{n+1},O_n
\overset{\iota_n}{\longrightarrow}
\hoA{\phi}/P_{n+1},O_n
\overset{\text{proj}}{\longrightarrow}
\oA(\phi)/D_{>n},P_{n+1},O_n
$$
where the first map is the map induced by $\iota_n$
which we also denote by $\iota_n$,
and the second map is the projection.
By Proposition \ref{prop.invariantbyP}
the equivalence class $[\cZ(L)]$
in the first set $\hoA{\coprod^l S^1}/L_{<2n},P_{n+1},O_n$ is
invariant under both of orientation change and Kirby move II.
Hence it is also true for the equivalence class
$[\iota_n \cZ(L)]$
in the third set $\oA(\phi)/D_{>n},P_{n+1},O_n$,
which is naturally isomorphic to the set $\cA(\phi)/D_{>n}$
by Lemma \ref{lem.Pvanishes}.
This completes the proof.
\end{pf}

Using the map $\hat\D$ we will show a proof of the following lemma
in the following section.

\begin{lem}\label{lem.invertible}
Let $U_+$ (resp. $U_-$) be the trivial knot with $+1$ (resp. $-1$) framing.
Then $[\iota_n(\cZ(U_\pm))]$ is invertible in $\cA(\phi)/D_{>n}$,
in which we define an algebra structure
such that the disjoint union of two chord diagrams is
the product of the diagrams.
\end{lem}

Using the above proposition and lemma,
we obtain the following theorem.

\begin{thm}
Let $L$ be an oriented framed link,
and $M$ the 3-manifold obtained by Dehn surgery on $S^3$ along $L$.
Then the equivalence class
$$
[\iota_n(\cZ(U_+)]^{-\s_+}
[\iota_n(\cZ(U_-)]^{-\s_-}
[\iota_n(\cZ(L))] \in \cA(\phi)/D_{>n}
$$
is a topological invariant of $M$
for any positive integer $n$,
where we denote by $\s_+$ (resp. $\s_-$)
the number of positive (resp. negative) eigenvalues
of the linking matrix of $L$.
\end{thm}

\begin{pf}
We obtain
invariance under both of
orientation change of $L$ and
Kirby move II
by Proposition \ref{prop.invariant}.
Note that we can apply Kirby move II in any way
though we prove the proposition
for Kirby move II preserving the orientation of $L$,
because we also showed the invariance
under orientation change.

We also obtain invariance under Kirby move I,
since the change of $[\iota_n \cZ(L)]$ under the move
cancels with the change of $\s_{\pm}$.
\end{pf}

\begin{defn}
We denote the above invariant by $\Omega_n(M)$.
\end{defn}

\section{A universal quantum invariant of 3-manifolds}

In this section we unify the series $\Omega_n(M)$
into a invariant $\Omega(M)$.
We show that
the series $\Omega_n(M)$ satisfies a property,
and that by the property
the series has the same information
(modulo the order of the first homology group)
 as $\Omega(M)$ has.
We further show that
the invariant $\Omega(M)$ satisfies a property
derived from the above property of the series,
and that there exists the logarithm of $\Omega(M)$
by the property;
we denote by the logarithm
a universal quantum invariant of 3-manifolds.

\subsection{A group-like property of the series $\Omega_n(M)$}

We denote by $\hD_{n_1,n_2}$
the map
$\cA(\phi)/D_{>n_1+n_2}\rightarrow
\cA(\phi)/D_{>n_1} \otimes \cA(\phi)/D_{>n_2}$
naturally induced by $\hD$;
note that this map is well defined
since the degree is preserved by $\hD$,
where we regard the sum of degrees
as the degree in the tensor product.
Using Theorem \ref{thm.bunshin} we obtain the following proposition.

\begin{prop}\label{prop.bunshin}
$$
\hD_{n_1,n_2}(\Omega_{n_1+n_2}(M)) = \Omega_{n_1}(M) \otimes \Omega_{n_2}(M).
$$
\end{prop}

\begin{pf}
Noting that $\hD$ is an algebra homomorphism in $\cA(\phi)$,
this proposition is a direct conclusion of
Theorem \ref{thm.bunshin} and the following lemma.
\end{pf}

\begin{lem}\label{lem.Dandiota}
Let $n$, $n_1$ and $n_2$ be
positive integers satisfying $n=n_1+n_2$.
Then the following diagram is commutative;
\begin{equation}\label{eq.Dandiota}
\begin{CD}
{\hoA{\coprod^l S^1}/L_{<2n},PO_n} @>{\iota_n}>>
\hoA{\phi}/PO_n \\
@V{\hD_{n_1,n_2}}VV
@VV{\hD}V \\
\hoA{\coprod^l S^1}/L_{<2n_1},PO_{n_1} \otimes
\hoA{\coprod^l S^1}/L_{<2n_2},PO_{n_2}
@>{\iota_{n_1} \otimes \iota_{n_2}}>>
\hoA{\phi}/PO_{n_1} \otimes \hoA{\phi}/PO_{n_2}
\end{CD}
\end{equation}
where we mean by $PO_k$ the equivalence relation
generated by $P_{k+1}$ and $O_k$.
\end{lem}

\begin{pf}
By Lemma \ref{lem.decreaselegs}
we can assume that
an element in $\hoA{\coprod^l S^1}/L_{<2n},PO_n$
is a disjoint union of
dashed trivalent graphs and
$l$ copies of a solid circle with $n$ dashed isolated chords;
in fact, the space consists of linear sums of such elements.
Since $\iota_{\star}$ is trivial for dashed trivalent graphs,
(\ref{eq.Dandiota}) is commutative for them.
Hence it is sufficient to show that
(\ref{eq.Dandiota}) is commutative for
a solid circle with $n$ dashed isolated chords;
we put it to be $D_n$.

The image of $D_n$ by the map $\iota_n$ is equal to
the chord diagram
obtained by attaching $n$ isolated chords to $T_{2n}^n$,
by the definition of $\iota_n$.
By the same argument in Figure \ref{fig.Pvanish},
$T_{2n}^n$ with one dashed isolated chords is equal to
$T_{2n-2}^{n-1}$ times a sum of
$2n-2$ and a dashed circle;
the sum is equal to $-2$ with the relation $O_n$.
Repeating this argument,
we can show that
the image is equal to
$(-2)(-4)\cdots (-2n) = (-2)^n n!$.
Therefore the clockwise image of $D_n$ in the diagram
is equal to $(-2)^n n!$.

The image of $D_n$ by the map $\hD$ is equal to
the sum of $\binom{n}{k} D_k \otimes D_{n-k}$
by the definition of $\hD$.
Note that it vanishes with $L_{<2n_1}$ or $L_{<2n_2}$
unless $k=n_1$.
Hence the image by $\hD_{n_1,n_2}$ is equal to
$\binom{n}{n_1} D_{n_1} \otimes D_{n_2}$.
By the same argument as above,
the image of it by the map $\iota_{n_1} \otimes \iota_{n_2}$
is equal to
$\binom{n}{n_1} (-2)^{n_1} n_1! (-2)^{n_2} n_2! = (-2)^n n!$
using $n=n_1+n_2$;
this coincides the above value.
\end{pf}

\begin{pf*}{Proof of Lemma \ref{lem.invertible}}
Applying Theorem \ref{thm.bunshin} to $\cZ(U_{\pm})$, $n-1$ times,
we have
the formula $\hD^{(n-1)}(\cZ(U_{\pm}))=(\cZ(U_{\pm}))^{\otimes n}$,
where we define the map
$\hD^{(k)} : \hcA(X) \to \hcA(X)^{\otimes (k+1)}$ by
$\hD^{(k)}= (\hD \otimes 1)\circ \hD^{(k-1)}$ recursively.
Further, putting $X=\phi$,
the map $\hD^{(k)}$ naturally induces the map
of $\cA(\phi)/D_{>k+1}$ to $(\cA(\phi)/D_{>1})^{\otimes (k+1)}$;
we denote it by $\hD^{(k)}_{1,1,\cdots,1}$.
Applying Lemma \ref{lem.Dandiota} to the above formula $n-1$ times,
we obtain the first equality of the following formula;
$$
\hD^{(n-1)}_{1,1,\cdots,1}([\iota_n(\cZ(U_{\pm}))])
=(\iota_1(\cZ(U_{\pm})))^{\otimes n}
= (\mp 1 + \frac{\theta}{16})^{\otimes n}
\in (\cA(\phi)/D_{>1})^{\otimes n}.
$$
Here the second equality is derived from
Lemma \ref{lem.iota1U} below.
Therefore the constant term of $[\iota_n(\cZ(U_{\pm}))]$
does not vanish,
which means that it is invertible.
\end{pf*}

We have the following lemma proved in \cite{LMMOII}.

\begin{lem}[\cite{LMMOII}]\label{lem.iota1U}
The following formula holds;
$$
[\iota_1(\cZ(U_{\pm}))]= \mp 1 + \frac{\theta}{16}
\in \cA(\phi)/D_{>1}
$$
where we mean by $\theta$
the dashed trivalent graph
consisting of three edges and two vertices,
as the Greek character $\theta$.
\end{lem}

\begin{rem}
In fact,
the above values of $\cZ(U_{\pm})$
are $-2$ times the values in \cite{LMMOII}.
The difference occurs from
the definition of the map $\iota_1$;
it was defined by
simply removing $\Theta$ components in \cite{LMMOII},
where we mean by $\Theta$
a solid circle with one dashed line
as the Greek letter $\Theta$.
On the other hand,
we define $\iota_1$ here
in the way how we replace
$\Theta$ with the trivial dashed circle,
which is equivalent to $-2$.
\end{rem}

\subsection{A power series invariant $\Omega(M)$}

We define a map $\e : \cA(\phi) \to \C$
to be the projection to the degree $0$ part
of a linear sum of chord diagrams;
recall that we regard the empty diagram as $1$
whose degree is $0$.
By the definitions of $\hD_{1,n-1}$ and $\e$,
we immediately obtain the following lemma.

\begin{lem}\label{lem.contraction}
The following formula holds;
$$
(\e \otimes 1)\circ \hD_{1,n-1} = p_{n,n-1}
$$
where we mean by $1$ the identity map and
we denote by $p_{n,n-1}$
the projection of $\cA(\phi)/D_{>n}$
to $\cA(\phi)/D_{>n-1}$.
\end{lem}

\begin{lem}\label{lem.unifyOmega}
Let $(\Omega_1,\Omega_2,\cdots)$ be a series of
$\Omega_n \in \cA(\phi)/D_{>n}$
which satisfies
\begin{equation}\label{eq.split}
\hD_{n_1,n_2}(\Omega_{n_1+n_2})=\Omega_{n_1} \otimes \Omega_{n_2}
\end{equation}
for any positive integers $n_1$ and $n_2$.
\begin{enumerate}
\item
Then the formula $\Omega_n^{(d)} = m^{n-d} \Omega_d^{(d)}$ holds
for $d<n$,
where we denote by $\al^{(d)}$ the degree $d$ part of $\al$
and we put $m$ to be $\Omega_1^{(0)}$.
\item
Further, for $\Omega \in \hat\cA(\phi)$
defined to be $1+ \sum_{n=1}^{\infty} \Omega_n^{(n)}$,
the formula $\hD(\Omega)=\Omega \otimes \Omega$ holds;
recall that we denote by $\hcA(\phi)$
the completion of $\cA(\phi)$
with respect to the degree,
that is, $\hcA(\phi)$
consists of infinite linear sums of chord diagrams.
We also denote by $\hD$
the natural extension of $\hD$ to $\hcA(\phi)$.
\end{enumerate}
\end{lem}

\begin{pf}
We apply the map in Lemma \ref{lem.contraction}
to $\Omega_n$.
Then we have the left hand side as;
$$
(\e \otimes 1)\circ \hD_{1,n-1} (\Omega_n)
= (\e \otimes 1) (\Omega_1 \otimes \Omega_{n-1})
= m \Omega_{n-1}.
$$
Hence we obtain (1).

We put $\Omega_0$ to be $1$,
then the series $(\Omega_0,\Omega_1,\cdots)$
satisfies (\ref{eq.split})
for any non-negative integers $n_1$ and $n_2$.
It is sufficient to show the required formula of (2)
in each degree $n$.
Hence we will show the formula
$\hD(\Omega_n^{(n)}) = \sum_{k_1+k_2=n} \Omega^{k_1} \otimes \Omega^{k_2}$.
Though this is a formula in $\hcA(\phi)\otimes\hcA(\phi)$,
it consists only of degree $n$ part.
Therefore it suffices to show it
for the image of it by the map $p_{\infty,n_1} \otimes p_{\infty,n_2}$
for each pair $n_1$ and $n_2$ with $n_1+n_2=n$,
where we denote by $p_{\infty,k}$
the projection of $\hcA(\phi)$ to $\cA(\phi)/D_{>k}$.
The image becomes
$\hD_{n_1,n_2}(\Omega_n^{(n)})=
\Omega_{n_1}^{(n_1)} \otimes \Omega_{n_2}^{(n_2)}$
which is a special case of (\ref{eq.split}).
This completes the proof.
\end{pf}

We have the following lemma by results in \cite{LMMO}
where only the invariant $\Omega_1(M)$ was discussed.

\begin{lem}\label{lem.degreezero}
The degree zero part of $\Omega_1(M)$ is equal to
$|H_1(M,\Z)|$ if $M$ is a rational homology 3-sphere,
$0$ otherwise.
Here we mean by $|\cdot|$ the order of the set.
\end{lem}

We know that we can apply Lemma \ref{lem.unifyOmega}
to our series of invariants $\Omega_n(M)$
by Proposition \ref{prop.bunshin},
to reduce the series to $\Omega \in \hA$ and a scalar $m$.
In our case the scalar $m$ becomes
either the order of the first homology group or zero
by Lemma \ref{lem.degreezero}.
Then we obtain the following definition and proposition
from Lemma \ref{lem.unifyOmega} (2).

\begin{defn}
We define a topological invariant $\Omega(M)$
of a 3-manifold $M$ by
$\Omega(M) = 1 + \sum_{n=1}^{\infty} \Omega_n(M)^{(n)} \in \hA$.
\end{defn}

\begin{prop}\label{prop.Omegasplit}
The invariant $\Omega(M)$ satisfies
$\hD(\Omega(M))=\Omega(M) \otimes \Omega(M)$.
\end{prop}

\subsection{Logarithm of $\Omega(M)$}

We denote by $\Ac$
the vector subspace of $\cA(\phi)$ spanned by
the set of connected non-empty dashed trivalent graphs
with oriented trivalent vertices
subject to AS and IHX relations.
We put $\hAc$ to be the completion of it
with respect to the degree,
which becomes a vector subspace of $\hA$.

It is well known
as a property of Hopf algebra,
see for example \cite{Hopfalgebra},
that a non-zero element $\Omega \in \hA$ is group-like
if and only if
there exists a primitive element $\omega \in \hA$ satisfying
$\Omega = \exp(\omega) = 1+\omega+(1/2)\omega^2+\cdots$,
where
we call $\al$ {\it group-like}
if it satisfies $\hD(\al)= \al \otimes \al$,
and call $\al$ {\it primitive}
if it satisfies $\hD(\al)=\al \otimes 1 + 1 \otimes \al$.
In our case $\al \in \hA$ is primitive
if and only if $\al$ belongs to $\hAc$.
Note that $\omega$ is uniquely determined for given $\Omega$,
since a primitive element always has a positive degree.

Since $\Omega(M)$ is group-like
by Proposition \ref{prop.Omegasplit},
we have the following definition.

\begin{defn}
We define $\omega(M) \in \hAc$
by $\exp(\omega(M))=\Omega(M)$.
We call $\omega(M)$ a universal quantum invariant.
\end{defn}

\section{Properties of the universal quantum invariant $\omega(M)$}

We will see some properties of $\omega(M)$ in this section.

\subsection{Formulas for connected sum and opposite orientation}

\begin{prop}
Let $M$ be the connected sum of
two closed 3-manifolds $M_1$ and $M_2$,
then $\omega(M)$ is given by
$$
\omega(M) = \sum_{d=1}^{\infty}
\left( m_2^{d} \omega(M_1)^{(d)} + m_1^{d} \omega(M_2)^{(d)} \right),
$$
where we put $m_i$ ($i=1,2$) to be
the order of $H_1(M_i,\Z)$
if $M_i$ is a rational homology 3-sphere,
$0$ otherwise.
\end{prop}

\begin{pf}
Let $L$ be the framed link
obtained by taking split union of
two framed links $L_1$ and $L_2$.
By definition of $\cZ(L)$
we have $\cZ(L)$
as the disjoint union of $\cZ(L_1)$ and $\cZ(L_2)$.

Note that,
if $M_1$ and $M_2$ are obtained by Dehn surgery along
$L_1$ and $L_2$ respectively, then
$M$ is obtained from $L$.
By definition of $\Omega_n(M)$ we have
$\Omega_n(M)=\Omega_n(M_1)\Omega_n(M_2)$,
recall that we define multiplication by disjoint union
of chord diagrams.
By Lemmas \ref{lem.unifyOmega} and \ref{lem.degreezero}
we have
\begin{align*}
\Omega_n(M)^{(n)}
&= \sum_{d_1+d_2=n} \Omega_n(M_1)^{(d_1)}\Omega_n(M_2)^{(d_2)} \\
&= \sum_{d_1+d_2=n} m_1^{d_2}\Omega_{d_1}(M_1)^{(d_1)}
                    m_2^{d_1}\Omega_{d_2}(M_1)^{(d_2)}.
\end{align*}
Hence we have
$$
\Omega(M) = \sum_{d_1,d_2=0}^{\infty}
m_1^{d_2} \Omega(M_1)^{(d_1)}  m_2^{d_1} \Omega(M_2)^{(d_2)}.
$$

If both of $M_1$ and $M_2$ are
rational homology 3-spheres,
the above formula implies that
$\sum_d \Omega(M_i)^{(d)}/m_i^d$
is multiplicative with respect to connected sum.
Hence $\sum_d \omega(M_i)^{(d)}/m_i^d$ is additive,
and we obtain the required formula.

If either of $M_1$ and $M_2$, say $M_1$,
is not a rational homology 3-sphere,
then we have $m_1=0$.
Hence
$\Omega(M)=\sum_d \Omega(M_1)^{(d)} m_2^d$
by the above formula putting $d_2=0$,
where we regard $0^0$ as $1$.
Therefore we obtain the required formula,
completing the proof.
\end{pf}

\begin{prop}
Let $-M$ be a 3-manifold $M$ with opposite orientation.
Then the following formula holds;
$$
\omega(-M) =
\sum_{d=1}^{\infty} (-1)^d \omega(M)^{(d)}.
$$
\end{prop}

\begin{pf}
We define a map $\hat S:\cA(X) \to \cA(X)$
by putting $\hat S(D)=(-1)^d D$ for a chord diagram $D$
where $d$ is the degree of $D$.
For the mirror image $\bar L$ of a framed link $L$,
we have $\cZ(\bar L) = \hat S(\cZ(L))$;
we can show the formula by checking it
for each elementary tangles.
Since $-M$ is obtained by Dehn surgery along $\bar L$,
we have
$\Omega_n(-M)=\hat S(\Omega_n(M))$,
$\Omega(-M)=\hat S(\Omega(M))$ and
$\omega(-M)=\hat S(\omega(M))$ by their definitions.
The last formula is the required one.
\end{pf}

\subsection{The first term in $\omega(M)$}

By results in \cite{LMMOII} we have the following proposition.

\begin{prop}
Let $M$ be a oriented closed 3-manifold.
Then the coefficient of dashed $\theta$-curve in $\omega(M)$ is equal to
$(-1)^{b_1(M)+1}3 \tilde\lambda (M)$
where $b_1(M)$ is the first Betti number of $M$ and
we denote $\tilde\lambda(M)$ twice Lescop's generalization \cite{Lescop}
of the Casson-Walker invariant (Walker's normalization,
which is twice Casson's normalization) $\lambda(M)$
\cite{AkbulutMcCarthy,Walker}
satisfying $\tilde\lambda(M)=|H_1(M;\Z)|\lambda(M)$
if $M$ is a rational homology 3-sphere.
\end{prop}

\begin{pf}
In \cite{LMMOII} it is shown that
the degree $1$ part of $\Omega_1(M)$
is equal to $(-1)^{b_1(M)+1}3 \tilde\lambda (M)$.
The proposition immediately follows from that.
\end{pf}

\begin{rem}
We can \lq\lq substitute'' a Lie algebra into dashed lines,
see \cite{BarNatanII}.
When we substitute $sl_N$, $so_N$ and $sp_N$
to dashed $\theta$-curve,
we obtain values $2N(N^2-1)$, $N(N-1)(N-2)/2$ and $2N(N+1)(2N+1)$
respectively.
When we substitute $sl_2$ into
the first term (i.e. dashed $\theta$-curve and its coefficient)
of $\omega(M)$,
we have the value $12$ times the coefficient.
Hence, if one expects the existence of
\lq\lq $G$ Casson invariant $\lambda^G(M)$''
and expects that it should recover
from the first term in $\omega(M)$
by substituting the Lie algebra of $G$ into dashed lines,
it must satisfy
\begin{align*}
\lambda^{SU(N)}(M) &= \frac{N(N^2-1)}{6} \lambda(M), \\
\lambda^{SO(N)}(M) &= \frac{N(N-1)(N-2)}{24} \lambda(M), \\
\lambda^{Sp(N)}(M) &= \frac{N(N+1)(2N+1)}{6} \lambda(M).
\end{align*}
In \cite{Ohtsuki}
we have the same formula for $\lambda^{SU(N)}(M)$ of each lens space $M$,
which is obtained expanding quantum $PSU(N)$ invariant of lens spaces
obtained in \cite{Takata}
into a power series in $q-1$.
\end{rem}

\end{document}